# Some infinite series involving the Riemann zeta function

Donal F. Connon

dconnon@btopenworld.com

16 May 2010

**Abstract**

This paper considers some infinite series involving the Riemann zeta function. Some examples are set out below

$$\sum_{n=0}^{\infty}\frac{1}{n+1}\sum_{k=0}^{n}\binom{n}{k}\frac{(-1)^k x^k}{(k+y)^s} = s\,\Phi(x,s+1,y) - \log x\,\Phi(x,s,y)$$

$$\sum_{n=0}^{\infty}(-1)^n\sum_{k=0}^{n}\binom{n}{k}\frac{2^k}{(k+1)^s} = \frac{1}{2}\left(1-2^{1-s}\right)\varsigma(s)$$

$$\sum_{n=0}^{\infty}\frac{1}{2^{n+1}}\sum_{k=0}^{n}\binom{n}{k}\left[\frac{(-1)^k}{k+y} + \frac{(-1)^k}{k+1-y}\right] = \pi\,\mathrm{cosec}(\pi y)$$

$$\sum_{n=1}^{\infty}\frac{2^{2n}\sin^{2n}u}{n^3}\binom{2n}{n}^{-1} = 4u^2\log(2\sin u) + 2\mathrm{Cl}_3(2u) + 4u\,\mathrm{Cl}_2(2u) - 2\varsigma(3)$$

$$\frac{7}{4}\varsigma(3) - \frac{1}{2}\pi G = \sum_{n=0}^{\infty}\frac{2^{4n}}{(2n+1)^3}\frac{[n!]^4}{[(2n)!]^2}$$

where $\Phi(x,s,y)$ is the Hurwitz-Lerch zeta function and $\mathrm{Cl}_n(t)$ are the Clausen functions.

**1. Some Hasse-type series**

It was shown in equation (4.4.85) in [21] (as recently corrected) that

(1.1) $\qquad \dfrac{y}{s-1}\sum_{n=0}^{\infty}\dfrac{1}{n+1}\sum_{k=0}^{n}\binom{n}{k}\dfrac{(-1)^k y^k}{(k+1)^{s-1}} = Li_s(y) - \dfrac{\log y}{s-1}Li_{s-1}(y)$

where $Li_s(y)$ is the polylogarithm function [39]

$$Li_s(y) = \sum_{n=1}^{\infty}\frac{y^n}{n^s}$$

and, with $y = 1$, this becomes the formula originally discovered by Hasse [32] in 1930

(1.2) $$\frac{1}{s-1}\sum_{n=0}^{\infty}\frac{1}{n+1}\sum_{k=0}^{n}\binom{n}{k}\frac{(-1)^k}{(k+1)^{s-1}} = Li_s(1) = \varsigma(s)$$

The gamma function is defined as

$$\Gamma(s) = \int_0^{\infty} t^{s-1}e^{-t}dt \quad , \operatorname{Re}(s) > 0$$

and, using the substitution $t = (k+y)u$, we obtain

(1.3) $$\frac{1}{(k+y)^s} = \frac{1}{\Gamma(s)}\int_0^{\infty} u^{s-1}e^{-(k+y)u}du$$

We now consider the finite sum set out below

(1.4) $$S_n(x, y) = \sum_{k=0}^{n}\binom{n}{k}\frac{x^k}{(k+y)^s}$$

and combine (1.3) and (1.4) to obtain

$$S_n(x, y) = \sum_{k=0}^{n}\binom{n}{k}\frac{x^k}{(k+y)^s} = \sum_{k=0}^{n}\binom{n}{k}x^k\frac{1}{\Gamma(s)}\int_0^{\infty}u^{s-1}e^{-(k+y)u}du$$

$$= \frac{1}{\Gamma(s)}\int_0^{\infty}u^{s-1}\left\{\sum_{k=0}^{n}\binom{n}{k}x^k e^{-(k+y)u}\right\}du$$

$$= \frac{1}{\Gamma(s)}\int_0^{\infty}u^{s-1}e^{-yu}\left\{\sum_{k=0}^{n}\binom{n}{k}\left[xe^{-u}\right]^k\right\}du$$

Using the binomial theorem this becomes

(1.5) $$S_n(x, y) = \frac{1}{\Gamma(s)}\int_0^{\infty}u^{s-1}e^{-yu}\left(1+xe^{-u}\right)^n du$$

Making the summation, we see that

$$\sum_{n=0}^{\infty}t^n\sum_{k=0}^{n}\binom{n}{k}\frac{x^k}{(k+y)^s} = \frac{1}{\Gamma(s)}\sum_{n=0}^{\infty}t^n\int_0^{\infty}u^{s-2}e^{-yu}\left(1+xe^{-u}\right)^n du$$

The geometric series gives us for $\left|t\left(1+xe^{-u}\right)\right|<1$

$$\sum_{n=0}^{\infty}\left[t\left(1+xe^{-u}\right)\right]^{n}=\frac{1}{1-t\left(1+xe^{-u}\right)}$$

and we then have

(1.6) $$\sum_{n=0}^{\infty}t^{n}\sum_{k=0}^{n}\binom{n}{k}\frac{x^{k}}{(k+y)^{s}}=\frac{1}{\Gamma(s)}\int_{0}^{\infty}\frac{u^{s-1}e^{-yu}}{1-t\left(1+xe^{-u}\right)}du$$

We now integrate (1.6) with respect to $t$

$$\sum_{n=0}^{\infty}\frac{w^{n+1}}{n+1}\sum_{k=0}^{n}\binom{n}{k}\frac{x^{k}}{(k+y)^{s}}=\frac{1}{\Gamma(s)}\int_{0}^{w}dt\int_{0}^{\infty}\frac{u^{s-1}e^{-yu}}{1-t(1+xe^{-u})}du$$

$$=-\frac{1}{\Gamma(s)}\int_{0}^{\infty}\frac{u^{s-1}e^{-yu}\log\left[1-w(1+xe^{-u})\right]}{1+xe^{-u}}du$$

and obtain

(1.7) $$\sum_{n=0}^{\infty}\frac{w^{n+1}}{n+1}\sum_{k=0}^{n}\binom{n}{k}\frac{x^{k}}{(k+y)^{s}}=-\frac{1}{\Gamma(s)}\int_{0}^{\infty}\frac{u^{s-1}e^{-(y-1)u}\log\left[1-w(1+xe^{-u})\right]}{e^{u}+x}du$$

When $w=1$ and $x\to -x$ we get

(1.8) $$\sum_{n=0}^{\infty}\frac{1}{n+1}\sum_{k=0}^{n}\binom{n}{k}\frac{(-1)^{k}x^{k}}{(k+y)^{s}}=-\frac{1}{\Gamma(s)}\int_{0}^{\infty}\frac{u^{s-1}e^{-(y-1)u}\log\left[xe^{-u}\right]}{e^{u}-x}du$$

$$=\frac{1}{\Gamma(s)}\int_{0}^{\infty}\frac{u^{s}e^{-(y-1)u}}{e^{u}-x}du-\frac{\log x}{\Gamma(s)}\int_{0}^{\infty}\frac{u^{s-1}e^{-(y-1)u}}{e^{u}-x}du$$

We see from (1.3) that

$$\frac{x^{k}}{(k+y)^{s}}=\frac{x^{k}}{\Gamma(s)}\int_{0}^{\infty}u^{s-1}e^{-(k+y)u}du$$

and we have the summation

$$\sum_{k=0}^{\infty} \frac{x^k}{(k+y)^s} = \frac{1}{\Gamma(s)} \int_0^{\infty} u^{s-1} e^{-yu} \sum_{k=0}^{\infty} x^k e^{-ku} \, du$$

$$= \frac{1}{\Gamma(s)} \int_0^{\infty} \frac{u^{s-1} e^{-yu}}{1 - xe^{-u}} \, du$$

We therefore obtain the well-known formula [43, p.121] for the Hurwitz-Lerch zeta function $\Phi(x, s, y)$

$$(1.9) \qquad \Phi(x, s, y) = \sum_{k=0}^{\infty} \frac{x^k}{(k+y)^s} = \frac{1}{\Gamma(s)} \int_0^{\infty} \frac{u^{s-1} e^{-(y-1)u}}{e^u - x} \, du$$

and with $y = 1$ we obtain

$$\Phi(x, s, 1) = \sum_{k=0}^{\infty} \frac{x^k}{(k+1)^s} = \frac{1}{x} \sum_{k=1}^{\infty} \frac{x^k}{k^s} = \frac{Li_s(x)}{x}$$

giving us [46, p.280]

$$(1.10) \qquad Li_s(x) = \frac{x}{\Gamma(s)} \int_0^{\infty} \frac{u^{s-1}}{e^u - x} \, du$$

Reference to (1.8) then shows that

$$(1.11) \qquad \sum_{n=0}^{\infty} \frac{1}{n+1} \sum_{k=0}^{n} \binom{n}{k} \frac{(-1)^k x^k}{(k+y)^s} = s \Phi(x, s+1, y) - \log x \, \Phi(x, s, y)$$

With $x = 1$ we obtain Hasse's formula (1.2)

$$\sum_{n=0}^{\infty} \frac{1}{n+1} \sum_{k=0}^{n} \binom{n}{k} \frac{(-1)^k}{(k+y)^s} = s \Phi(1, s+1, y) = s\varsigma(s+1, y)$$

With $y = 1$ equation (1.11) becomes

$$\sum_{n=0}^{\infty} \frac{1}{n+1} \sum_{k=0}^{n} \binom{n}{k} \frac{(-1)^k x^k}{(k+1)^s} = s \Phi(x, s+1, 1) - \log x \, \Phi(x, s, 1)$$

or equivalently

(1.12) $$\sum_{n=0}^{\infty}\frac{1}{n+1}\sum_{k=0}^{n}\binom{n}{k}\frac{(-1)^k x^k}{(k+1)^s} = \frac{s}{x}Li_{s+1}(x) - \frac{1}{x}\log x\, Li_s(x)$$

which corresponds with (1.1). A closed form expression may be obtained for example with $x = 1/2$ and $s = 2$.

We see from (1.6) that

$$\sum_{n=0}^{\infty} t^n \sum_{k=0}^{n}\binom{n}{k}\frac{x^k}{(k+y)^s} = \frac{1}{(1-t)\Gamma(s)}\int_0^{\infty}\frac{u^{s-1}e^{-yu}}{1-e^{-u}tx/(1-t)}du$$

$$= \frac{1}{(1-t)\Gamma(s)}\int_0^{\infty}\frac{u^{s-1}e^{-(y-1)u}}{e^u - tx/(1-t)}du$$

and hence referring to (1.9) we have

(1.13) $$\sum_{n=0}^{\infty} t^n \sum_{k=0}^{n}\binom{n}{k}\frac{x^k}{(k+y)^s} = \frac{1}{1-t}\Phi\left(\frac{tx}{1-t}, s, y\right)$$

With $t = 1/2$ we obtain

(1.13.1) $$\sum_{n=0}^{\infty}\frac{1}{2^{n+1}}\sum_{k=0}^{n}\binom{n}{k}\frac{x^k}{(k+y)^s} = \Phi(x, s, y)$$

which is clearly an example of Euler's transformation of series [33, p.244].

Letting $w = \dfrac{tx}{1-t}$ (1.13) may be represented by

(1.14) $$\sum_{n=0}^{\infty}\left(\frac{w}{w+x}\right)^n \sum_{k=0}^{n}\binom{n}{k}\frac{x^k}{(k+y)^s} = \frac{w+x}{x}\Phi(w, s, y)$$

and with $x = -1$ this becomes

(1.15) $$\sum_{n=0}^{\infty}\left(\frac{-w}{1-w}\right)^n \sum_{k=0}^{n}\binom{n}{k}\frac{(-1)^k}{(k+y)^s} = (1-w)\Phi(w, s, y)$$

as previously reported by Guillera and Sondow [30i].

Letting $w = 1$ and $x = -1$ in (1.14) results in another identity due to Hasse

(1.16) $$\sum_{n=0}^{\infty}\frac{1}{2^{n+1}}\sum_{k=0}^{n}\binom{n}{k}\frac{(-1)^k}{(k+y)^s} = \Phi(-1,s,y) = \varsigma_a(s,y)$$

where $\varsigma_a(s,y)$ is the alternating Hurwitz-Lerch zeta function

$$\varsigma_a(s,y) = \sum_{k=0}^{\infty}\frac{(-1)^k}{(k+y)^s}$$

With $y=1$ we have

(1.17) $$\sum_{n=0}^{\infty}\frac{1}{2^{n+1}}\sum_{k=0}^{n}\binom{n}{k}\frac{(-1)^k}{(k+1)^s} = \varsigma_a(s,1) = \varsigma_a(s)$$

Letting $w = x = 1$ in (1.14) results in

(1.18) $$\sum_{n=0}^{\infty}\frac{1}{2^{n+1}}\sum_{k=0}^{n}\binom{n}{k}\frac{1}{(k+y)^s} = \Phi(1,s,y) = \varsigma(s,y)$$

Letting $s=1$ and $y=1/2$ in (1.16) results in

$$\sum_{n=0}^{\infty}\frac{1}{2^{n}}\sum_{k=0}^{n}\binom{n}{k}\frac{(-1)^k}{2k+1} = \varsigma_a(1,1/2)$$

We have from [15, p.523]

$$\varsigma_a(1,1/2) = \frac{\pi}{2}$$

and we therefore obtain

$$\sum_{n=0}^{\infty}\frac{1}{2^{n}}\sum_{k=0}^{n}\binom{n}{k}\frac{(-1)^k}{2k+1} = \frac{\pi}{2}$$

which was also determined in a different manner in equation (8.11e) in [22].

Since [15, p.523]

$$\varsigma_a(1,y) + \varsigma_a(1,1-y) = \pi\operatorname{cosec}(\pi y)$$

we have for $0 < y < 1$

$$\sum_{n=0}^{\infty} \frac{1}{2^{n+1}} \sum_{k=0}^{n} \binom{n}{k} \left[ \frac{(-1)^k}{k+y} + \frac{(-1)^k}{k+1-y} \right] = \pi \csc(\pi y)$$

□

We now divide (1.7) by $w$ and integrate to obtain

$$\sum_{n=0}^{\infty} \frac{v^{n+1}}{(n+1)^2} \sum_{k=0}^{n} \binom{n}{k} \frac{x^k}{(k+y)^s} = -\frac{1}{\Gamma(s)} \int_0^v \frac{\log\left[1 - w(1 + xe^{-u})\right]}{w} dw \int_0^{\infty} \frac{u^{s-2} e^{-(y-1)u}}{e^u + x} du$$

Since

$$\int \frac{\log(1-ax)}{x} dx = -Li_2(ax)$$

we have

(1.19) $$\sum_{n=0}^{\infty} \frac{v^{n+1}}{(n+1)^2} \sum_{k=0}^{n} \binom{n}{k} \frac{x^k}{(k+y)^s} = \frac{1}{\Gamma(s)} \int_0^{\infty} \frac{u^{s-1} e^{-(y-1)u} Li_2\left[v(1+xe^{-u})\right]}{e^u + x} du$$

Further operations of the same kind will result in

(1.20) $$\sum_{n=0}^{\infty} \frac{v^{n+1}}{(n+1)^r} \sum_{k=0}^{n} \binom{n}{k} \frac{x^k}{(k+y)^s} = \frac{1}{\Gamma(s)} \int_0^{\infty} \frac{u^{s-1} e^{-(y-1)u} Li_r\left[v(1+xe^{-u})\right]}{e^u + x} du$$

□

With $t = -1$ in (1.6) we obtain

$$\sum_{n=0}^{\infty} (-1)^n \sum_{k=0}^{n} \binom{n}{k} \frac{x^k}{(k+y)^s} = \frac{1}{\Gamma(s)} \int_0^{\infty} \frac{u^{s-1} e^{-yu}}{2 + xe^{-u}} du$$

and letting $x \to -x$ we get

$$\sum_{n=0}^{\infty} (-1)^n \sum_{k=0}^{n} \binom{n}{k} \frac{(-1)^k x^k}{(k+y)^s} = \frac{1}{\Gamma(s)} \int_0^{\infty} \frac{u^{s-1} e^{-yu}}{2 - xe^{-u}} du$$

With $x = 2$ we have

$$\sum_{n=0}^{\infty} (-1)^n \sum_{k=0}^{n} \binom{n}{k} \frac{(-1)^k 2^k}{(k+y)^s} = \frac{1}{2\Gamma(s)} \int_0^{\infty} \frac{u^{s-1} e^{-yu}}{1 - e^{-u}} du = \frac{1}{2\Gamma(s)} \int_0^{\infty} \frac{u^{s-1} e^{-(y-1)u}}{e^u - 1} du$$

Hence we obtain using (1.9)

$$(1.21) \qquad \sum_{n=0}^{\infty}(-1)^n \sum_{k=0}^{n}\binom{n}{k}\frac{(-1)^k 2^k}{(k+y)^s} = \frac{1}{2}\Phi(1,s,y)$$

and with $s = n+1$ and $y = 1$ this becomes

$$(1.22) \qquad \sum_{n=0}^{\infty}(-1)^n \sum_{k=0}^{n}\binom{n}{k}\frac{(-1)^k 2^k}{(k+1)^{n+1}} = \frac{1}{2}\varsigma(n+1)$$

This identity, in the case where $n$ is a positive integer, was recently reported by Alzer and Koumandos [6] where it was derived in a very different manner.

With $x = -2$ we have

$$\sum_{n=0}^{\infty}(-1)^n \sum_{k=0}^{n}\binom{n}{k}\frac{2^k}{(k+y)^s} = \frac{1}{2\Gamma(s)}\int_0^{\infty}\frac{u^{s-2}e^{-yu}}{e^u+1}du$$

$$= \frac{1}{2}\Phi(-1,s,y)$$

we see that with $y = 1$

$$= \frac{1}{2}\Phi(-1,s,1) = \frac{1}{2}\varsigma_a(s)$$

Hence we have

$$(1.23) \qquad \sum_{n=0}^{\infty}(-1)^n \sum_{k=0}^{n}\binom{n}{k}\frac{2^k}{(k+1)^s} = \frac{1}{2}\left(1-2^{2-s}\right)\varsigma(s)$$

With $t = x = \frac{1}{2}$ in (1.6) we obtain

$$\sum_{n=0}^{\infty}\frac{1}{2^n}\sum_{k=0}^{n}\binom{n}{k}\frac{1}{2^k(k+y)^s} = \frac{2}{\Gamma(s)}\int_0^{\infty}\frac{u^{s-2}e^{-yu}}{e^u-1}du$$

$$= 2\Phi(1,s,y) = 2\varsigma(s,y)$$

and therefore we have with $y = 1$

(1.24) $$\varsigma(s) = \sum_{n=0}^{\infty} \frac{1}{2^{n+1}} \sum_{k=0}^{n} \binom{n}{k} \frac{1}{2^k (k+1)^s}$$

With $t = -1$, $x = 1/2$ and $y = 1$ in (1.6) we have

(1.25) $$\sum_{n=0}^{\infty} (-1)^n \sum_{k=0}^{n} \binom{n}{k} \frac{1}{2^k (k+1)^s} = \frac{1}{2\Gamma(s)} \int_0^{\infty} \frac{u^{s-1}}{e^u + \frac{1}{4}} du = -2 Li_s(-1/4)$$

## 2. Some series involving the central binomial numbers

From, for example, Knopp's book [33, p.271] we have the well-known Maclaurin expansion

(2.0) $$\left(\sin^{-1} y\right)^2 = \frac{1}{2} \sum_{n=1}^{\infty} \frac{[n!]^2}{n^2 (2n)!} (2y)^{2n} \quad , \quad -1 \leq y \leq 1$$

and therefore upon letting $x = \sin^{-1} y$ we get

(2.1) $$x^2 = \frac{1}{2} \sum_{n=1}^{\infty} \frac{[n!]^2}{n^2 (2n)!} 2^{2n} \sin^{2n} x \quad , \quad -\pi/2 \leq x \leq \pi/2$$

which was known by Euler [34].

Using Bürmann's theorem, it is an exercise in Whittaker & Watson [46, p.130] to prove

$$x^2 = \sin^2 x + \left(\frac{2}{3}\right) \frac{1}{2} \sin^4 x + \left(\frac{2.4}{3.5}\right) \frac{1}{3} \sin^6 x + ...$$

$$= \sum_{n=1}^{\infty} A_n \sin^{2n} x$$

where $$A_n = \frac{(2n-2)!!}{(2n-1)!!} \frac{1}{n}$$

Using the definitions of the double factorials

$$(2n)!! = 2.4...(2n) = 2^n n! \qquad (2n+1)!! = 1.3.5...(2n+1) = \frac{(2n+1)!}{2^n n!}$$

we see this is the same as the above identity.

We now multiply equation (2.1) by $\cot x$ and integrate to obtain

$$\int_0^t x^2 \cot x \, dx = \frac{1}{2} \sum_{n=1}^{\infty} \frac{[n!]^2}{n^2 (2n)!} 2^{2n} \int_0^t \sin^{2n-1} x \cos x \, dx$$

which, for $-\pi/2 \leq x \leq \pi/2$, results in [12, p.234]

(2.2) $$\int_0^t x^2 \cot x \, dx = \frac{1}{4} \sum_{n=1}^{\infty} \frac{2^{2n} \sin^{2n} t}{n^3} \binom{2n}{n}^{-1}$$

where $\binom{2n}{n} = \frac{(2n)!}{[n!]^2}$ is the central binomial coefficient.

This integral also appears in [12, p.234].

Using integration by parts we have

$$\int_0^t x^2 \cot x \, dx = x^2 \log \sin x \Big|_0^t - 2 \int_0^t x \log \sin x \, dx$$

Since $x^2 \log \sin x = x^2 \log \frac{\sin x}{x} + x^2 \log x$ we see that $\lim_{x \to 0} x^2 \log \sin x = 0$ and thus

$$\int_0^t x^2 \cot x \, dx = t^2 \log \sin t - 2 \int_0^t x \log \sin x \, dx$$

Hence we get for $0 \leq t \leq \pi/2$

(2.3) $$\int_0^t x \log \sin x \, dx = \frac{1}{2} t^2 \log \sin t - \frac{1}{8} \sum_{n=1}^{\infty} \frac{2^{2n} \sin^{2n} t}{n^3} \binom{2n}{n}^{-1}$$

According to Ayoub [8, p.1084], Euler returned to the zeta function, for what appears to be the last time, in 1772 in a paper entitled "Exercitationes Analyticae" [26]. Notwithstanding that by this time Euler had been blind for six years, through what Ayoub describes as "a striking and elaborate scheme", he was able to prove that

(2.4) $$\int_0^{\pi/2} x \log \sin x \, dx = \frac{7}{16} \varsigma(3) - \frac{\pi^2}{8} \log 2$$

A very elementary proof of (2.4) is given in equation (6.20a) in [22] where we used the basic identity

(2.5) $$\int_a^b p(x)\cot x\,dx = 2\sum_{n=1}^{\infty}\int_a^b p(x)\sin 2nx\,dx$$

which, as shown in [22], is valid for a wide class of suitably behaved functions. Specifically we require that $p(x)$ is a twice continuously differentiable function. It should be noted that in the above formula we require either (i) both $\sin(x/2)$ and $\cos(x/2)$ have no zero in $[a,b]$ or (ii) if either $\sin(a/2)$ or $\cos(a/2)$ is equal to zero then $p(a)$ must also be zero. Condition (i) is equivalent to the requirement that $\sin x$ has no zero in $[a,b]$.

For example, letting $p(x) = x^2$ in (2.5) we have for $-\pi < t < \pi$

$$\int_0^t x^2 \cot x\,dx = 2\sum_{n=1}^{\infty}\int_0^t x^2 \sin 2nx\,dx$$

$$= \frac{1}{2}\sum_{n=1}^{\infty}\frac{\cos 2nx}{n^3} - x^2\sum_{n=1}^{\infty}\frac{\cos 2nx}{n} + x\sum_{n=1}^{\infty}\frac{\sin 2nx}{n^2}\bigg|_0^t$$

This gives us

(2.6) $$\int_0^t x^2 \cot x\,dx = \frac{1}{2}\sum_{n=1}^{\infty}\frac{\cos 2nt}{n^3} - t^2\sum_{n=1}^{\infty}\frac{\cos 2nt}{n} + t\sum_{n=1}^{\infty}\frac{\sin 2nt}{n^2} - \frac{1}{2}\varsigma(3)$$

and with $t = \pi/2$ we obtain

$$\int_0^{\pi/2} x^2 \cot x\,dx = \frac{1}{2}\sum_{n=1}^{\infty}\frac{(-1)^n}{n^3} - \frac{\pi^2}{4}\sum_{n=1}^{\infty}\frac{(-1)^n}{n} - \frac{1}{2}\varsigma(3)$$

Hence using the alternating zeta function

$$\varsigma_a(s) = \sum_{n=1}^{\infty}\frac{(-1)^{n+1}}{n^s} = (1-2^{1-s})\varsigma(s)$$

we see that

(2.7) $$\int_0^{\pi/2} x^2 \cot x\,dx = -\frac{7}{8}\varsigma(3) + \frac{\pi^2}{4}\log 2$$

Another elementary evaluation of this integral has recently been provided by Fujii and Suzuki [28].

Using the substitutions $y = \sin x$ and $z = \tan x$ we may also note that

$$\int_0^t x^n \cot x \, dx = \int_0^{\sin t} \frac{\left(\sin^{-1} y\right)^n}{y} dy$$

$$2\int_0^t \frac{x^n}{\sin 2x} dx = \int_0^{\tan t} \frac{\left(\tan^{-1} z\right)^n}{z} dz$$

We also see from (2.2) and (2.6) that for $-\pi/2 \leq t \leq \pi/2$

(2.8) $\quad \dfrac{1}{4}\sum_{n=1}^{\infty}\dfrac{2^{2n}\sin^{2n}t}{n^3}\binom{2n}{n}^{-1} = \dfrac{1}{2}\sum_{n=1}^{\infty}\dfrac{\cos 2nt}{n^3} - t^2\sum_{n=1}^{\infty}\dfrac{\cos 2nt}{n} + t\sum_{n=1}^{\infty}\dfrac{\sin 2nt}{n^2} - \dfrac{1}{2}\varsigma(3)$

and with $t = \pi/2$ we immediately see that

(2.9) $\quad \varsigma(3) = \dfrac{2}{7}\pi^2 \log 2 - \dfrac{2}{7}\sum_{n=1}^{\infty}\dfrac{2^{2n}}{n^3}\binom{2n}{n}^{-1}$

which is in agreement with Sherman's compendium of formulae [42].

We have from [12, p.198]

(2.10) $\quad B(n, 1/2) = \dfrac{2^{2n}}{n}\binom{2n}{n}^{-1}$

where the beta function is defined by

$$B(u, v) = \int_0^1 (1-t)^{u-1} t^{v-1} dt$$

This then gives us

$$\dfrac{2^{2n}}{n}\binom{2n}{n}^{-1} = \int_0^1 \dfrac{(1-t)^n}{(1-t)\sqrt{t}} dt$$

and we have the summation involving the dilogarithm function

$$\sum_{n=1}^{\infty} \dfrac{2^{2n}}{n^3}\binom{2n}{n}^{-1} = \int_0^1 \dfrac{Li_2(1-t)}{(1-t)\sqrt{t}} dt$$

It may be noted that the Wolfram Mathematica Online Integrator evaluates this integral in terms of polylogarithms of order 2 and 3 and hence this will not provide us with any new information regarding the nature of $\varsigma(3)$. It should be easier to obtain a simpler version of the output by using integration by parts, noting that

$$\int \frac{Li_2(1-t)}{\sqrt{t}} dt = 2\int Li_2(1-u^2) du$$

and

$$\int Li_2(1-u^2) du = u Li_2(1-u^2) - Li_2(1-u) - Li_2(-u) - \log u \log(1+u) - 2u \log u - 2u$$

We now attempt to use (2.10) in (2.2)

$$\int_0^t x^2 \cot x \, dx = \frac{1}{4}\sum_{n=1}^{\infty} \frac{2^{2n} \sin^{2n} t}{n^3} \binom{2n}{n}^{-1}$$

and this gives us

$$\sum_{n=1}^{\infty} \frac{2^{2n} \sin^{2n} t}{n^3} \binom{2n}{n}^{-1} = \int_0^1 \frac{Li_2\left[(1-t)\sin^2 t\right]}{(1-t)\sqrt{t}} dt$$

Following an idea used by Batir [9] we now substitute an integral for the dilogarithm function; we have

(2.10.1) $\quad Li_n(z) = \frac{(-1)^{n-1}}{(n-1)!}\int_0^1 \frac{z\log^{n-1} u}{1-zu} du \quad\quad Li_2(z) = -\int_0^1 \frac{z\log u}{1-zu} du$

and thus

$$\int_0^1 \frac{Li_2\left[(1-t)\sin^2 t\right]}{(1-t)\sqrt{t}} dt = -\int_0^1\int_0^1 \frac{\sin^2 t \log u}{\sqrt{t}\left[1-u(1-t)\sin^2 t\right]} dt\, du$$

However the Wolfram Mathematica Online Integrator cannot evaluate this integral.

The integral (2.10.1) was also used by Batir [9] to derive (2.73) using his complex double integral

$$\sum_{n=1}^{\infty} \frac{x^{2n}[n!]^4}{(2n)^k[(2n)!]^2} = \frac{(-1)^{k-1}}{(k-3)!}\int_0^{\pi/2} \frac{1}{\sin y}\int_0^{\arcsin(x\sin y/4)} \phi \log^{k-3}[4\sin\phi/x\sin y]d\phi dy$$

Since $\sin\left(\dfrac{\pi}{6}\right) = \dfrac{1}{2}$ we see from (2.2) that

$$(2.11) \qquad \sum_{n=1}^{\infty} \frac{1}{n^3}\binom{2n}{n}^{-1} = -8\int_0^{\pi/6} x\log[2\sin x]\,dx$$

which I first came across in van der Poorten's 1979 paper "Some wonderful formulae…an introduction to Polylogarithms" [45].

In passing, we note from [12, p.122] that

$$\sum_{n=1}^{\infty} \frac{(2t)^{2n}}{n^3}\binom{2n}{n}^{-1} = 4\int_0^t \frac{\left(\sin^{-1} x\right)^2}{x}\,dx = 2t^2\,_4F_3\left[\{1,1,1,1\},\left\{\frac{3}{2},2,2\right\},t^2\right]$$

in terms of the hypergeometric functions.

□

We now multiply (2.1) by $x\cot x$ and integrate

$$\int_0^t x^3 \cot x\,dx = \frac{1}{2}\sum_{n=1}^{\infty} \frac{[n!]^2}{n^2(2n)!} 2^{2n} \int_0^t x\sin^{2n-1} x \cos x\,dx$$

Integration by parts yields

$$\int_0^t x\sin^{2n-1} x \cos x\,dx = \frac{t}{2n}\sin^{2n} t - \frac{1}{2n}\int_0^t \sin^{2n} x\,dx$$

and therefore we have

$$(2.11.1) \qquad \int_0^t x^3 \cot x\,dx = \frac{1}{4}t\sum_{n=1}^{\infty} \frac{2^{2n}\sin^{2n} t}{n^3}\binom{2n}{n}^{-1} - \frac{1}{4}\sum_{n=1}^{\infty} \frac{2^{2n}}{n^3}\binom{2n}{n}^{-1}\int_0^t \sin^{2n} x\,dx$$

With $t = \pi/2$ we obtain

$$\int_0^{\pi/2} x^3 \cot x\,dx = \frac{\pi}{8}\sum_{n=1}^{\infty} \frac{2^{2n}}{n^3}\binom{2n}{n}^{-1} - \frac{1}{4}\sum_{n=1}^{\infty} \frac{2^{2n}}{n^3}\binom{2n}{n}^{-1}\int_0^{\pi/2} \sin^{2n} x\,dx$$

and using [12, p.195]

(2.11.2) $$\int_0^{\pi/2} \sin^{2n} x \, dx = \frac{(2n)!}{2^{2n}(n!)^2} \frac{\pi}{2} = \frac{1}{2^{2n}} \binom{2n}{n} \frac{\pi}{2}$$

we have

$$\int_0^{\pi/2} x^3 \cot x \, dx = \frac{\pi}{8} \sum_{n=1}^{\infty} \frac{2^{2n}}{n^3} \binom{2n}{n}^{-1} - \frac{\pi}{8} \sum_{n=1}^{\infty} \frac{1}{n^3}$$

Using (2.9) we obtain

(2.12) $$\int_0^{\pi/2} x^3 \cot x \, dx = \frac{\pi^3}{8} \log 2 - \frac{9\pi}{16} \varsigma(3)$$

As another example, letting $p(x) = x^3$ in (2.5) we have for $-\pi < t < \pi$

$$\int_0^t x^3 \cot x \, dx = 2 \sum_{n=1}^{\infty} \int_0^t x^3 \sin 2nx \, dx$$

$$= \frac{3}{2} x^2 \sum_{n=1}^{\infty} \frac{\sin 2nx}{n^2} - \frac{3}{4} \sum_{n=1}^{\infty} \frac{\sin 2nx}{n^4} - x^3 \sum_{n=1}^{\infty} \frac{\cos 2nx}{n} + \frac{3}{2} x \sum_{n=1}^{\infty} \frac{\cos 2nx}{n^3} \bigg|_0^t$$

$$= \frac{3}{2} t^2 \sum_{n=1}^{\infty} \frac{\sin 2nt}{n^2} - \frac{3}{4} \sum_{n=1}^{\infty} \frac{\sin 2nt}{n^4} - t^3 \sum_{n=1}^{\infty} \frac{\cos 2nt}{n} + \frac{3}{2} t \sum_{n=1}^{\infty} \frac{\cos 2nt}{n^3}$$

Using (2.11.1) together with (2.56) we see that

$$\int_0^t x^3 \cot x \, dx = \frac{1}{4} t \sum_{n=1}^{\infty} \frac{2^{2n} \sin^{2n} t}{n^3} \binom{2n}{n}^{-1} - \frac{1}{4} \sum_{n=1}^{\infty} \frac{2^{2n}}{n^3} \binom{2n}{n}^{-1} \frac{1}{2^{2n}} \left[ \binom{2n}{n} t + \sum_{j=1}^{n} \frac{(-1)^j}{j} \binom{2n}{n-j} \sin 2jt \right]$$

$$= \frac{1}{4} t \sum_{n=1}^{\infty} \frac{2^{2n} \sin^{2n} t}{n^3} \binom{2n}{n}^{-1} - \frac{t}{4} \sum_{n=1}^{\infty} \frac{1}{n^3} - \frac{1}{4} \sum_{n=1}^{\infty} \frac{1}{n^3} \binom{2n}{n}^{-1} \sum_{j=1}^{n} \frac{(-1)^j}{j} \binom{2n}{n-j} \sin 2jt$$

This gives us

$$\frac{1}{4} t \sum_{n=1}^{\infty} \frac{2^{2n} \sin^{2n} t}{n^3} \binom{2n}{n}^{-1} - \frac{t}{4} \varsigma(3) - \frac{1}{4} \sum_{n=1}^{\infty} \frac{1}{n^3} \binom{2n}{n}^{-1} \sum_{j=1}^{n} \frac{(-1)^j}{j} \binom{2n}{n-j} \sin 2jt$$

$$= \frac{3}{2}t^2 \sum_{n=1}^{\infty} \frac{\sin 2nt}{n^2} - \frac{3}{4}\sum_{n=1}^{\infty} \frac{\sin 2nt}{n^4} - t^3 \sum_{n=1}^{\infty} \frac{\cos 2nt}{n} + \frac{3}{2}t \sum_{n=1}^{\infty} \frac{\cos 2nt}{n^3}$$

and upon substituting (2.8) we have for $-\pi/2 \le t \le \pi/2$

(2.12.1)

$$\sum_{n=1}^{\infty} \frac{1}{n^3} \binom{2n}{n}^{-1} \sum_{j=1}^{n} \frac{(-1)^j}{j} \binom{2n}{n-j} \sin 2jt = 3\sum_{n=1}^{\infty} \frac{\sin 2nt}{n^4} - 4t \sum_{n=1}^{\infty} \frac{\cos 2nt}{n^3} - 2t^2 \sum_{n=1}^{\infty} \frac{\sin 2nt}{n^2} - 3t\varsigma(3)$$

□

It was shown in equation (4.3.168a) of [19] that for $0 \le t < 1$

(2.13)

$$\int_0^t \pi x^2 \cot \pi x\, dx = -[\varsigma'(-2,t) + \varsigma'(-2,1-t)] + 2t[\varsigma'(-1,t) - \varsigma'(-1,1-t)] + t^2 \log(2\sin \pi t) - \frac{\varsigma(3)}{2\pi^2}$$

where $\varsigma'(s,t) = \frac{\partial}{\partial s}\varsigma(s,t)$ is a derivative of the Hurwitz zeta function.

Using (2.2) we see that for $0 \le t \le 1/2$

(2.13.1) $\qquad \int_0^t \pi x^2 \cot \pi x\, dx = \frac{1}{\pi^2}\int_0^{\pi t} u^2 \cot u\, du = \frac{1}{4\pi^2} \sum_{n=1}^{\infty} \frac{2^{2n} \sin^{2n} \pi t}{n^3}\binom{2n}{n}^{-1}$

and hence we see from (2.13) and (2.13.1) that for $0 \le t \le 1/2$
(2.14)

$$\frac{1}{4\pi^2}\sum_{n=1}^{\infty} \frac{2^{2n} \sin^{2n} \pi t}{n^3}\binom{2n}{n}^{-1} = -[\varsigma'(-2,t) + \varsigma'(-2,1-t)] + 2t[\varsigma'(-1,t) - \varsigma'(-1,1-t)] + t^2 \log(2\sin \pi t) - \frac{\varsigma(3)}{2\pi^2}$$

I initially thought of letting $t \to 1-t$ in (2.14) with the hope that

$$\frac{1}{4\pi^2}\sum_{n=1}^{\infty} \frac{2^{2n} \sin^{2n} \pi t}{n^3}\binom{2n}{n}^{-1} = -[\varsigma'(-2,t) + \varsigma'(-2,1-t)] - 2(1-t)[\varsigma'(-1,t) - \varsigma'(-1,1-t)]$$

$$+(1-t)^2 \log(2\sin \pi t) - \frac{\varsigma(3)}{2\pi^2}$$

and subtraction of these two equations appeared to indicate that for

$$\varsigma'(-1,t)-\varsigma'(-1,1-t)=\frac{1}{2}(1-2t)\log(2\sin \pi t)$$

We note from (2.22) below that

$$\varsigma'\left(-1,\frac{1}{4}\right)-\varsigma'\left(-1,\frac{3}{4}\right)=\frac{1}{2\pi}\mathrm{Cl}_2(\pi/2)=\frac{G}{2\pi}$$

and apparently this would lead to a closed-form expression for Catalan's constant. Unfortunately, the analysis is incorrect because the restriction that $0\leq t\leq 1/2$ prevents us from letting $t\to 1-t$ in (2.14).

We may also express (2.13) in terms of the Barnes multiple gamma functions $\Gamma_n(x)$ defined in [43, p.24]. Adamchik [5] has shown that for $\mathrm{Re}(x)>0$

(2.15) $$\varsigma'(-n,x)-\varsigma'(-n)=(-1)^n\sum_{k=0}^{n}k!Q_{k,n}(x)\log\Gamma_{k+1}(x)$$

where $Q_{k,n}(x)$ are polynomials defined by

$$Q_{k,n}(x)=\sum_{j=k}^{n}(1-x)^{n-j}\binom{n}{j}\begin{Bmatrix}j\\k\end{Bmatrix}$$

and $\begin{Bmatrix}j\\k\end{Bmatrix}$ are the Stirling subset numbers defined by

$$\begin{Bmatrix}j\\k\end{Bmatrix}=k\begin{Bmatrix}n-1\\k\end{Bmatrix}+\begin{Bmatrix}n-1\\k-1\end{Bmatrix},\quad \begin{Bmatrix}n\\0\end{Bmatrix}=\begin{cases}1, & n=0\\0, & n\neq 0\end{cases}$$

Particular cases of (2.15) are

(2.16) $$\varsigma'(-1,x)-\varsigma'(-1)=x\log\Gamma(x)-\log G(x+1)$$

$$\varsigma'(-2,x)-\varsigma'(-2)=2\log\Gamma_3(x)+(3-2x)\log G(x)-(1-x)^2\log\Gamma(x)$$

and letting $x\to 1-x$ we see that

$$\varsigma'(-1,1-x)-\varsigma'(-1)=(1-x)\log\Gamma(1-x)-\log\Gamma(1-x)-\log G(1-x)$$

$$=-x\log\Gamma(1-x)-\log G(1-x)$$

$$\varsigma'(-2,1-x)-\varsigma'(-2)=2\log\Gamma_3(1-x)+(1+2x)\log G(1-x)-x^2\log\Gamma(1-x)$$

Hence we have

(2.17) $$\varsigma'(-1,x)-\varsigma'(-1,1-x)=x\log[\Gamma(x)\Gamma(1-x)]+\log\frac{G(1-x)}{G(1+x)}$$

which we have previously derived in equation (4.3.162) of [19].

Letting $t=\frac{1}{2}$ in (2.14) we obtain

$$\frac{1}{4\pi^2}\sum_{n=1}^{\infty}\frac{2^{2n}}{n^3}\binom{2n}{n}^{-1}=-2\varsigma'\left(-2,\frac{1}{2}\right)+\frac{1}{4}\log 2-\frac{\varsigma(3)}{2\pi^2}$$

and using (2.9) we obtain the known result

(2.18) $$\varsigma'\left(-2,\frac{1}{2}\right)=\frac{3\varsigma(3)}{16\pi^2}$$

which is also derived, for example, in equation (4.3.168d) of [19].

When $t=\frac{1}{4}$ in (2.14) we get

(2.19)
$$\frac{1}{4\pi^2}\sum_{n=1}^{\infty}\frac{2^{n}}{n^3}\binom{2n}{n}^{-1}=-\left[\varsigma'\left(-2,\frac{1}{4}\right)+\varsigma'\left(-2,\frac{3}{4}\right)\right]+\frac{1}{2}\left[\varsigma'\left(-1,\frac{1}{4}\right)-\varsigma'\left(-1,\frac{3}{4}\right)\right]+\frac{1}{32}\log 2-\frac{\varsigma(3)}{2\pi^2}$$

We now refer to Adamchik's result [2]

(2.20) $$\varsigma'(-2n-1,t)-\varsigma'(-2n-1,1-t)=\frac{(2n+1)!}{(2\pi)^{2n+1}}\text{Cl}_{2n+2}(2\pi t)$$

(2.21) $$\varsigma'(-2n,t)+\varsigma'(-2n,1-t)=(-1)^n\frac{(2n)!}{(2\pi)^{2n}}\text{Cl}_{2n+1}(2\pi t)$$

where $\text{Cl}_n(t)$ is the Clausen function defined by

(2.21.1) $$\text{Cl}_{2n+1}(t)=\sum_{k=1}^{\infty}\frac{\cos kt}{k^{2n+1}} \qquad \text{Cl}_{2n}(t)=\sum_{k=1}^{\infty}\frac{\sin kt}{k^{2n}}$$

This is also derived in equation (4.3.167) of [19].

We then see that (as noted by Adamchik [2])

(2.22) $$\varsigma'\left(-1,\frac{1}{4}\right)-\varsigma'\left(-1,\frac{3}{4}\right)=\frac{1}{2\pi}Cl_2(\pi/2)=\frac{G}{2\pi}$$

where G is Catalan's constant.

Using the formula [18, equation (5.11)] (where I have corrected a misprint)

(2.22.1) $$Cl_{2n+1}(\pi/2)=\frac{1-2^{2n}}{2^{4n+1}}\varsigma(2n+1)$$

this then gives us for $t=1/4$

$$\varsigma'\left(-2n,\frac{1}{4}\right)+\varsigma'\left(-2n,\frac{3}{4}\right)=(-1)^n\frac{(2n)!}{(2\pi)^{2n}}Cl_{2n+1}(\pi/2)=(-1)^n\frac{1-2^{2n}}{2^{4n+1}}\frac{(2n)!}{(2\pi)^{2n}}\varsigma(2n+1)$$

and in particular we have for $n=1$

$$\varsigma'\left(-2,\frac{1}{4}\right)+\varsigma'\left(-2,\frac{3}{4}\right)=\frac{3\varsigma(3)}{64\pi^2}$$

We therefore have

(2.23) $$\sum_{n=1}^{\infty}\frac{2^n}{n^3}\binom{2n}{n}^{-1}=-\frac{35\varsigma(3)}{16}+\frac{\pi^2}{8}\log 2+\pi G$$

in agreement with equation (3.167) of Sherman's paper [42].

We note that this also concurs with equation (6.69o) of [22] and also [18] where, by entirely different methods, it is shown that

(2.24) $$\int_0^{\pi/4} x^2 \cot x\, dx = -\frac{35}{64}\varsigma(3)+\frac{\pi^2}{32}\log 2+\frac{\pi G}{4}$$

With $t=\frac{1}{6}$ in (2.14) we obtain

(2.25) $$\frac{1}{4\pi^2}\sum_{n=1}^{\infty}\frac{1}{n^3}\binom{2n}{n}^{-1}=-\left[\varsigma'\left(-2,\frac{1}{6}\right)-\varsigma'\left(-2,\frac{5}{6}\right)\right]+\frac{1}{3}\left[\varsigma'\left(-1,\frac{1}{6}\right)-\varsigma'\left(-1,\frac{5}{6}\right)\right]-\frac{\varsigma(3)}{2\pi^2}$$

We see from (2.21) that

$$\varsigma'\left(-2,\frac{1}{6}\right)+\varsigma'\left(-2,\frac{5}{6}\right)=-\frac{1}{2\pi^2}Cl_3\left(\frac{\pi}{3}\right)$$

and from Lewin's monograph [39, p.198] we have

$$Cl_{2n+1}\left(\frac{\pi}{3}\right)=\frac{1}{2}(1-2^{-2n})(1-3^{-2n})\varsigma(2n+1)$$

which gives us

$$Cl_3\left(\frac{\pi}{3}\right)=\frac{1}{4}\varsigma(3)$$

This results in

(2.26) $$\varsigma'\left(-2,\frac{1}{6}\right)+\varsigma'\left(-2,\frac{5}{6}\right)=-\frac{\varsigma(3)}{8\pi^2}$$

We also have from (2.20)

$$\varsigma'\left(-1,\frac{1}{6}\right)-\varsigma'\left(-1,\frac{5}{6}\right)=\frac{1}{2\pi}Cl_2\left(\frac{\pi}{3}\right)$$

and from (2.25) we obtain

(2.27) $$\frac{1}{4\pi^2}\sum_{n=1}^{\infty}\frac{1}{n^3}\binom{2n}{n}^{-1}=-\left[\varsigma'\left(-2,\frac{1}{6}\right)-\varsigma'\left(-2,\frac{5}{6}\right)\right]+\frac{1}{6\pi}Cl_2\left(\frac{\pi}{3}\right)-\frac{\varsigma(3)}{2\pi^2}$$

According we have two simultaneous equations involving $\varsigma'\left(-2,\frac{1}{6}\right)$ and $\varsigma'\left(-2,\frac{5}{6}\right)$ (unfortunately (2.27) contains two other unknown constants; in this regard see (2.33)).

We have the well-known Hurwitz's formula for the Fourier expansion of the Riemann zeta function $\varsigma(s,t)$ as reported in Titchmarsh's treatise [43b, p.37]

(2.28) $$\varsigma(s,t)=2\Gamma(1-s)\left[\sin\left(\frac{\pi s}{2}\right)\sum_{n=1}^{\infty}\frac{\cos 2n\pi t}{(2\pi n)^{1-s}}+\cos\left(\frac{\pi s}{2}\right)\sum_{n=1}^{\infty}\frac{\sin 2n\pi t}{(2\pi n)^{1-s}}\right]$$

where $\operatorname{Re}(s) < 0$ and $0 < t \leq 1$. In 2000, Boudjelkha [12b] showed that this formula also applies in the region $\operatorname{Re}(s) < 1$. It may be noted that when $t = 1$ this reduces to Riemann's functional equation for $\varsigma(s)$. Letting $s \to 1-s$ we may write this as

$$(2.29) \quad \varsigma(1-s,t) = 2\Gamma(s)\left[\cos\left(\frac{\pi s}{2}\right)\sum_{n=1}^{\infty}\frac{\cos 2n\pi t}{(2\pi n)^s} + \sin\left(\frac{\pi s}{2}\right)\sum_{n=1}^{\infty}\frac{\sin 2n\pi t}{(2\pi n)^s}\right]$$

Using (2.29) it is easily shown that

$$(2.30) \quad 2\varsigma'(-1,t) - B_2(t)(1-\gamma-\log(2\pi)) = -4\sum_{n=1}^{\infty}\frac{\log n \cos 2n\pi t}{(2\pi n)^2} + \frac{1}{2\pi}\operatorname{Cl}_2(2\pi t)$$

$$(2.31) \quad \varsigma'(-2,t) - B_3(t)\left[\frac{1}{2} - \frac{1}{3}\gamma - \frac{1}{2}\log(2\pi)\right] = -4\sum_{n=1}^{\infty}\frac{\log n \sin 2n\pi t}{(2\pi n)^3} - \frac{1}{(2\pi)^2}\operatorname{Cl}_3(2\pi t)$$

With $t = 1/6$ in (2.31) we obtain

$$\varsigma'\left(-2,\frac{1}{6}\right) - \frac{5}{216}\left[\frac{1}{2} - \frac{1}{3}\gamma - \frac{1}{2}\log(2\pi)\right] = -4\sum_{n=1}^{\infty}\frac{\log n \sin(n\pi/3)}{(2\pi n)^3} - \frac{\varsigma(3)}{16\pi^2}$$

and this indicates the complexities involved in determining a closed form expression for $\varsigma'\left(-2,\frac{1}{6}\right)$.

Differentiating (2.39) gives us

$$-\sum_{n=1}^{\infty}\frac{\log n \sin(n\pi/3)}{n^s} = \sqrt{3}\left\{\frac{3^{-s}-1}{2}\varsigma'(s) - \frac{3^{-s}\log 3}{2}\varsigma(s)\right\}$$

$$+\sqrt{3}\left\{6^{-s}\left[\varsigma'\left(s,\frac{1}{6}\right) + \varsigma'\left(s,\frac{1}{3}\right)\right] - 6^{-s}\log 6\left[\varsigma\left(s,\frac{1}{6}\right) + \varsigma\left(s,\frac{1}{3}\right)\right]\right\}$$

so that with $s = 3$ we have

$$-\sum_{n=1}^{\infty}\frac{\log n \sin(n\pi/3)}{n^3} = \sqrt{3}\left\{\frac{3^{-3}-1}{2}\varsigma'(3) - \frac{3^{-3}\log 3}{2}\varsigma(3)\right\}$$

$$+\sqrt{3}\left\{6^{-3}\left[\varsigma'\left(3,\frac{1}{6}\right) + \varsigma'\left(3,\frac{1}{3}\right)\right] - 6^{-3}\log 6\left[\varsigma\left(3,\frac{1}{6}\right) + \varsigma\left(3,\frac{1}{3}\right)\right]\right\}$$

With $t = 1/4$ in (2.30) we obtain

$$2\varsigma'\left(-1,\frac{1}{4}\right) + \frac{1}{48}[1-\gamma-\log(2\pi)] = -4\sum_{n=1}^{\infty}\frac{\log n \cos(n\pi/2)}{(2\pi n)^2} + \frac{1}{2\pi}Cl_2\left(\frac{\pi}{2}\right)$$

$$2\varsigma'\left(-1,\frac{1}{4}\right) + \frac{1}{48}[1-\gamma-\log(2\pi)] = -4\sum_{n=1}^{\infty}(-1)^n\frac{\log 2n}{(4\pi n)^2} + \frac{G}{2\pi}$$

$$= \frac{\log 2}{4\pi^2}\sum_{n=1}^{\infty}\frac{(-1)^{n+1}}{n^2} + \frac{1}{4\pi^2}\sum_{n=1}^{\infty}(-1)^{n+1}\frac{\log n}{n^2} + \frac{G}{2\pi}$$

$$= \frac{\log 2}{4\pi^2}\varsigma_a(2) - \frac{1}{4\pi^2}\varsigma'_a(2) + \frac{G}{2\pi}$$

Since $\varsigma_a(2) = \frac{1}{2}\varsigma(2)$ and $\varsigma'_a(2) = \frac{1}{2}\varsigma(2)\log 2 + \frac{1}{2}\varsigma'(2)$ we have

$$2\varsigma'\left(-1,\frac{1}{4}\right) + \frac{1}{48}[1-\gamma-\log(2\pi)] = -\frac{1}{8\pi^2}\varsigma'(2) + \frac{G}{2\pi}$$

It is easily found from the functional equation for the Riemann zeta function that

$$\varsigma'(-1) = \frac{1}{12}(1-\gamma-\log(2\pi)) + \frac{1}{2\pi^2}\varsigma'(2)$$

and hence we have (as originally determined by Adamchik [2])

(2.32) $$\varsigma'\left(-1,\frac{1}{4}\right) = \frac{G}{4\pi} - \frac{1}{8}\varsigma'(-1)$$

□

Ghusayni [29] showed in 1998 that

(2.33) $$\varsigma(3) = \frac{\pi}{2}\sum_{n=1}^{\infty}\frac{\sin(n\pi/3)}{n^2} - \frac{3}{4}\sum_{n=1}^{\infty}\frac{1}{n^3}\binom{2n}{n}^{-1}$$

Later in 2000, Ghusayni [30], having noted an earlier paper [31], reported that

(2.34) $$\sum_{n=1}^{\infty}\frac{\sin(n\pi/3)}{n^2} = \frac{\sqrt{3}}{2}\left\{\frac{1}{1^2} + \frac{1}{2^2} - \frac{1}{4^2} - \frac{1}{5^2} + \frac{1}{7^2} + \frac{1}{8^2} - \ldots\right\}$$

$$= \frac{\sqrt{3}}{2}\left\{\frac{1}{3}\psi^{(1)}\left(\frac{1}{3}\right) - \frac{2}{9}\pi^2\right\}$$

and Ghusayni concluded that

(2.35) $$\varsigma(3) = -\frac{\sqrt{3}}{18}\pi^3 + \frac{3\sqrt{3}}{4}\pi\sum_{n=1}^{\infty}\frac{1}{(3n-2)^2} - \frac{3}{4}\sum_{n=1}^{\infty}\frac{1}{n^3}\binom{2n}{n}^{-1}$$

It is shown in [30] that

$$\frac{1}{1^2} + \frac{1}{2^2} - \frac{1}{4^2} - \frac{1}{5^2} + \frac{1}{7^2} + \frac{1}{8^2} - \ldots = -\frac{2}{27}\pi^2 - 2\int_0^1 \frac{\log x}{1+x^3}dx$$

*Mathematica* evaluates this integral as

$$324\int_0^1 \frac{\log x}{1+x^3}dx = -8\pi^2 - 6\text{PolyGamma}\left[1,\frac{1}{6}\right] + 3\text{PolyGamma}\left[1,\frac{5}{6}\right]$$

$$-3\text{PolyGamma}\left[1,\frac{1}{3}\right] + 6\text{PolyGamma}\left[1,\frac{2}{3}\right]$$

The Mathworld website for the central binomial coefficient reports the following formula which was experimentally obtained by Plouffe [41] (and in fact determined analytically in 1985 by Zucker [49])

(2.36) $$\sum_{n=1}^{\infty}\frac{1}{n^3}\binom{2n}{n}^{-1} = \frac{1}{18}\pi\sqrt{3}\left[\psi^{(1)}\left(\frac{1}{3}\right) - \psi^{(1)}\left(\frac{2}{3}\right)\right] - \frac{4}{3}\varsigma(3)$$

We have the reflection formula [43, p.14]

$$\psi(1-x) - \psi(x) = \pi\cot\pi x$$

and differentiation gives us

$$\psi^{(k)}(1-x) + (-1)^{k+1}\psi^{(k)}(x) = (-1)^k \pi\frac{d^k}{dx^k}\cot\pi x$$

and therefore we see that

$$\psi^{(1)}\left(\frac{1}{3}\right) + \psi^{(1)}\left(\frac{2}{3}\right) = \pi^2/\sin^2\left(\frac{\pi}{3}\right) = \frac{4}{3}\pi^2$$

We then obtain

(2.37) $$\sum_{n=1}^{\infty} \frac{1}{n^3} \binom{2n}{n}^{-1} = \frac{1}{9} \pi \sqrt{3} \, \psi^{(1)}\left(\frac{1}{3}\right) - \frac{4}{3} \varsigma(3) - \frac{2}{27} \pi^3 \sqrt{3}$$

We have

$$\psi^{(k)}\left(\frac{p}{q}\right) = (-1)^{k+1} k! \, \varsigma\left(k+1, \frac{p}{q}\right)$$

and therefore we see that

$$\psi'\left(\frac{1}{3}\right) = \varsigma\left(2, \frac{1}{3}\right) = 9 \sum_{n=0}^{\infty} \frac{1}{(3n+1)^2}$$

or, equivalently, changing the order of summation

$$\psi'\left(\frac{1}{3}\right) = 9 \sum_{n=1}^{\infty} \frac{1}{(3n-2)^2}$$

Hence we see that (2.36) is equivalent to Ghusayni's result (2.35)

$$\varsigma(3) = -\frac{\sqrt{3}}{18} \pi^3 + \frac{\sqrt{3}}{12} \pi \psi^{(1)}\left(\frac{1}{3}\right) - \frac{3}{4} \sum_{n=1}^{\infty} \frac{1}{n^3} \binom{2n}{n}^{-1}$$

□

We now refer back to (2.8) which is valid for $-1/2 \leq t \leq 1/2$

$$\frac{1}{4} \sum_{n=1}^{\infty} \frac{2^{2n} \sin^{2n} \pi t}{n^3} \binom{2n}{n}^{-1} = \frac{1}{2} \sum_{n=1}^{\infty} \frac{\cos 2n\pi t}{n^3} - \pi^2 t^2 \sum_{n=1}^{\infty} \frac{\cos 2n\pi t}{n} + \pi t \sum_{n=1}^{\infty} \frac{\sin 2n\pi t}{n^2} - \frac{1}{2} \varsigma(3)$$

so that with $t = 1/6$ we have

$$\frac{1}{4} \sum_{n=1}^{\infty} \frac{1}{n^3} \binom{2n}{n}^{-1} = \frac{1}{2} \sum_{n=1}^{\infty} \frac{\cos(n\pi/3)}{n^3} - \frac{1}{36} \pi^2 \sum_{n=1}^{\infty} \frac{\cos(n\pi/3)}{n} + \frac{1}{6} \pi \sum_{n=1}^{\infty} \frac{\sin(n\pi/3)}{n^2} - \frac{1}{2} \varsigma(3)$$

Lewin [39] and Srivastava and Tsumura [43, p.293] reported for $\text{Re}(s) > 1$

(2.38) $$\sum_{n=1}^{\infty} \frac{\cos(n\pi/3)}{n^s} = \frac{1}{2}(6^{1-s} - 3^{1-s} - 2^{1-s} + 1)\varsigma(s)$$

$$(2.39) \quad \sum_{n=1}^{\infty} \frac{\sin(n\pi/3)}{n^s} = \sqrt{3}\left\{\frac{3^{-s}-1}{2}\varsigma(s) + 6^{-s}\left[\varsigma\left(s,\frac{1}{6}\right) + \varsigma\left(s,\frac{1}{3}\right)\right]\right\}$$

Hence we have with $s = 3$ and $s = 2$ respectively

$$(2.40) \quad \sum_{n=1}^{\infty} \frac{\cos(n\pi/3)}{n^3} = \frac{1}{3}\varsigma(3)$$

$$(2.41) \quad \sum_{n=1}^{\infty} \frac{\sin(n\pi/3)}{n^2} = \sqrt{3}\left\{\frac{3^{-2}-1}{2}\varsigma(2) + 6^{-2}\left[\varsigma\left(2,\frac{1}{6}\right) + \varsigma\left(2,\frac{1}{3}\right)\right]\right\}$$

We have the Fourier series [44, p.148] for $0 < t < 1$

$$\sum_{n=1}^{\infty} \frac{\cos 2n\pi t}{n} = -\log(2\sin \pi t)$$

Using (2.41) we then obtain Ghusayni's result

$$\varsigma(3) = \frac{\pi}{2}\sqrt{3}\left\{\frac{3^{-2}-1}{2}\varsigma(2) + 6^{-2}\left[\varsigma\left(2,\frac{1}{6}\right) + \varsigma\left(2,\frac{1}{3}\right)\right]\right\} - \frac{3}{4}\sum_{n=1}^{\infty} \frac{1}{n^3}\binom{2n}{n}^{-1}$$

Srivastava and Tsumura [43, p.293] have also reported for $\text{Re}(s) > 1$

$$(2.42) \quad \sum_{n=1}^{\infty} \frac{\cos(2n\pi/3)}{n^s} = \frac{1}{2}(3^{1-s}-1)\varsigma(s)$$

$$(2.43) \quad \sum_{n=1}^{\infty} \frac{\sin(2n\pi/3)}{n^s} = \sqrt{3}\left\{\frac{3^{-s}-1}{2}\varsigma(s) + 3^{-s}\varsigma\left(s,\frac{1}{3}\right)\right\}$$

$$(2.44) \quad \sum_{n=1}^{\infty} \frac{\cos(n\pi/2)}{n^s} = 2^{-s}(2^{1-s}-1)\varsigma(s)$$

$$(2.45) \quad \sum_{n=1}^{\infty} \frac{\sin(n\pi/2)}{n^s} = (2^{-s}-1)\varsigma(s) + 2^{1-2s}\varsigma\left(s,\frac{1}{4}\right)$$

and the relevant values may be easily inserted in (2.8).

Reference should also be made to the recent paper "Certain series related to the triple sine function" by Koyama and Kurokawa [34].

□

More generally we have from (2.14), (2.20) and (2.21)

$$(2.46) \quad \sum_{n=1}^{\infty} \frac{2^{2n} \sin^{2n} \pi t}{n^3} \binom{2n}{n}^{-1} = 4\pi^2 t^2 \log(2 \sin \pi t) + 2\text{Cl}_3(2\pi t) + 4\pi t\, \text{Cl}_2(2\pi t) - 2\varsigma(3)$$

or equivalently

$$(2.47) \quad \sum_{n=1}^{\infty} \frac{2^{2n} \sin^{2n} u}{n^3} \binom{2n}{n}^{-1} = 4u^2 \log(2 \sin u) + 2\text{Cl}_3(2u) + 4u\, \text{Cl}_2(2u) - 2\varsigma(3)$$

Having recently seen a paper by Bradley et al. [14], I noted that this is equivalent to the identity previously discovered by Zucker [49] in 1985 (and I thank John Zucker for subsequently sending me a reprint of his original paper).

With $u = \pi/6$ and $u = \pi/4$ we obtain respectively

$$(2.48) \quad \sum_{n=1}^{\infty} \frac{1}{n^3} \binom{2n}{n}^{-1} = 2\text{Cl}_3(\pi/3) + \frac{2\pi}{3} \text{Cl}_2(\pi/3) - 2\varsigma(3)$$

$$(2.49) \quad \sum_{n=1}^{\infty} \frac{2^n}{n^3} \binom{2n}{n}^{-1} = \frac{\pi^2}{8} \log 2 + 2\text{Cl}_3(\pi/2) + \pi\, \text{Cl}_2(\pi/2) - 2\varsigma(3)$$

Then using (2.22.1) we see that

$$\text{Cl}_3(\pi/2) = -\frac{3}{32} \varsigma(3)$$

and hence from (2.49) we obtain another derivation of (2.23)

$$\sum_{n=1}^{\infty} \frac{2^n}{n^3} \binom{2n}{n}^{-1} = -\frac{35\varsigma(3)}{16} + \frac{\pi^2}{8} \log 2 + \pi G$$

We also have a number of formulae involving $\text{Cl}_2\left(\frac{p\pi}{q}\right)$ in Browkin's paper in Lewin's survey [39b, p.244], including for example

$$\text{Cl}_2\left(\frac{\pi}{6}\right) + \text{Cl}_2\left(\frac{5\pi}{6}\right) = \frac{4}{3} G = \frac{4}{3} \text{Cl}_2\left(\frac{\pi}{2}\right)$$

It is easily seen from the definition of the Clausen function (2.21.1) that

$$\text{Cl}_{2n}(\pi) = \text{Cl}_{2n}(2\pi) = 0$$

$$\text{Cl}_{2n+1}(\pi) = (2^{-2n} - 1)\varsigma(2n+1) = -\varsigma_a(2n+1)$$

$$\text{Cl}_{2n+1}(2\pi) = \varsigma(2n+1)$$

We also have

$$\text{Cl}_2(\pi/2) = G = -\text{Cl}_2(3\pi/2)$$

$$\frac{1}{2}\text{Cl}_2(2x) = \text{Cl}_2(x) - \text{Cl}_2(\pi - x)$$

which implies that

$$\text{Cl}_2\left(\frac{2\pi}{3}\right) = \frac{2}{3}\text{Cl}_2\left(\frac{\pi}{3}\right)$$

The Clausen function may be expressed in closed form in at least three other cases and from Lewin's book [39, p.198] we have

$$\text{Cl}_{2n+1}\left(\frac{\pi}{2}\right) = -2^{-2n-1}(1 - 2^{-2n})\varsigma(2n+1)$$

$$\text{Cl}_{2n+1}\left(\frac{\pi}{3}\right) = \frac{1}{2}(1 - 2^{-2n})(1 - 3^{-2n})\varsigma(2n+1)$$

$$\text{Cl}_{2n+1}\left(\frac{2\pi}{3}\right) = -\frac{1}{2}(1 - 3^{-2n})\varsigma(2n+1)$$

For example, we see from the definition that

$$\text{Cl}_{2n+1}\left(\frac{\pi}{2}\right) = -\frac{1}{2^{2n+1}} + \frac{1}{4^{2n+1}} - \frac{1}{6^{2n+1}} + \ldots$$

$$= -\frac{1}{2^{2n+1}}\left[\frac{1}{1^{2n+1}} - \frac{1}{2^{2n+1}} + \ldots\right] = -\frac{1}{2^{2n+1}}\varsigma_a(2n+1)$$

We then see $\text{Cl}_3(\pi/3) = \varsigma(3)/3$ and from (2.48) that

$$\sum_{n=1}^{\infty}\frac{1}{n^3}\binom{2n}{n}^{-1} = \frac{2\pi}{3}\text{Cl}_2(\pi/3) - \frac{4}{3}\varsigma(3)$$

and comparing this with (2.37) we obtain

(2.50) $\quad 2\sqrt{3}\,\mathrm{Cl}_2(\pi/3) = \psi^{(1)}\!\left(\dfrac{1}{3}\right) - \dfrac{2}{3}\pi^2$

as previously derived by Fettis [26a].

More generally we have [39a, p.358]

(2.51) $\quad \psi^{(1)}\!\left(\dfrac{p}{q}\right) = \dfrac{\pi^2}{2}\operatorname{cosec}^2\!\left(\dfrac{\pi p}{q}\right) + 2q \displaystyle\sum_{m=1}^{\lfloor (q-1)/2 \rfloor} \sin\!\left(\dfrac{2m\pi p}{q}\right)\mathrm{Cl}_2\!\left(\dfrac{2m\pi p}{q}\right)$

and a particular case of this is (2.50). Another example is

(2.52) $\quad \mathrm{Cl}_2(\pi/6) = \dfrac{1}{24}\!\left[\sqrt{3}\,\psi^{(1)}\!\left(\dfrac{1}{3}\right) + 16G - 2\pi^2/\sqrt{3}\right]$

Using PSLQ, Bailey et al. [8a] discovered experimentally that

(2.53) $\quad 6\,\mathrm{Cl}_2(\alpha) - 6\,\mathrm{Cl}_2(2\alpha) + 2\,\mathrm{Cl}_2(3\alpha) = 7\,\mathrm{Cl}_2\!\left(\dfrac{2\pi}{7}\right) + 7\,\mathrm{Cl}_2\!\left(\dfrac{4\pi}{7}\right) - 7\,\mathrm{Cl}_2\!\left(\dfrac{6\pi}{7}\right)$

where $\alpha = 2\tan^{-1}\sqrt{7}$. It appears that there must be a connection with (2.52) in the case where $q = 7$. See also Coffey's recent paper [18a].

□

Integrating (2.46) results in

(2.54) $\quad \displaystyle\sum_{n=1}^{\infty} \dfrac{2^{2n}}{n^3}\binom{2n}{n}^{-1} \int_0^x \sin^{2n} u\, du$

$$= \dfrac{4}{3}x^3 \log 2 + 4\int_0^x u^2 \log \sin u\, du + 2\,\mathrm{Cl}_4(2x) - 2x\,\mathrm{Cl}_3(2x) - 2x\,\varsigma(3)$$

In the above we used the fact that $\displaystyle\int_0^x u \sin 2ku\, du = \dfrac{\sin 2kx}{4k^2} - \dfrac{x\cos 2kx}{2k}$ and this results in

$$\int_0^x u\,\mathrm{Cl}_2(2u)\,du = \dfrac{1}{4}\mathrm{Cl}_4(2x) - \dfrac{1}{2}x\,\mathrm{Cl}_3(2x).$$

When $x = \pi/2$ we get

$$\sum_{n=1}^{\infty} \frac{2^{2n}}{n^3} \binom{2n}{n}^{-1} \int_0^{\pi/2} \sin^{2n} u \, du = \sum_{n=1}^{\infty} \frac{2^{2n}}{n^3} \binom{2n}{n}^{-1} \frac{\pi}{2} \frac{1}{2^{2n}} \binom{2n}{n} = \frac{\pi}{2} \sum_{n=1}^{\infty} \frac{1}{n^3} = \frac{\pi}{2} \varsigma(3)$$

$$= \frac{\pi^3}{6} \log 2 + 4 \int_0^{\pi/2} u^2 \log \sin u \, du - \pi \, \mathrm{Cl}_3(\pi) - \pi \varsigma(3)$$

We see that

$$\mathrm{Cl}_3(\pi) = \sum_{k=1}^{\infty} \frac{(-1)^k}{k^3} = -\varsigma_a(3) = -\frac{3}{4} \varsigma(3)$$

which gives us

$$\frac{\pi}{2} \varsigma(3) = \frac{\pi^3}{6} \log 2 + 4 \int_0^{\pi/2} u^2 \log \sin u \, du - \frac{\pi}{4} \varsigma(3)$$

and we therefore obtain a new derivation of Euler's integral

(2.55) $$\int_0^{\pi/2} u^2 \log \sin u \, du = -\frac{\pi^3}{24} \log 2 + \frac{3\pi}{16} \varsigma(3)$$

□

Wiener [47] has shown that

(2.56) $$\int_0^x \sin^{2n} u \, du = \frac{1}{2^{2n}} \left[ \binom{2n}{n} x + \sum_{j=1}^{n} \frac{(-1)^j}{j} \binom{2n}{n-j} \sin 2jx \right]$$

and we therefore obtain from (2.54)

(2.57) $$\sum_{n=1}^{\infty} \frac{1}{n^3} \binom{2n}{n}^{-1} \sum_{j=1}^{n} \frac{(-1)^j}{j} \binom{2n}{n-j} \sin 2jx$$

$$= \frac{4}{3} x^3 \log 2 + 4 \int_0^x u^2 \log \sin u \, du + 2\mathrm{Cl}_4(2x) - 2x \mathrm{Cl}_3(2x) - 3x \varsigma(3)$$

Using integration by parts we have

$$\int_0^x u^3 \cot u\, du = u^3 \log \sin u \Big|_0^x - 3\int_0^x u^2 \log \sin u\, du$$

$$= x^3 \log \sin x - 3\int_0^x u^2 \log \sin u\, du$$

Using (2.5) we have

$$\int_0^x u^3 \cot u\, du = 2\sum_{n=1}^\infty \int_0^x u^3 \sin 2nu\, du$$

$$= \frac{3}{4}u^2 \sum_{n=1}^\infty \frac{\sin 2nu}{n^2} - \frac{3}{8}\sum_{n=1}^\infty \frac{\sin 2nu}{n^4} - \frac{1}{2}u^3 \sum_{n=1}^\infty \frac{\cos 2nu}{n} + \frac{3}{4}u\sum_{n=1}^\infty \frac{\cos 2nu}{n^3}\bigg|_0^x$$

$$= \frac{3}{4}x^2 \sum_{n=1}^\infty \frac{\sin 2nx}{n^2} - \frac{3}{8}\sum_{n=1}^\infty \frac{\sin 2nx}{n^4} - \frac{1}{2}x^3 \sum_{n=1}^\infty \frac{\cos 2nx}{n} + \frac{3}{4}x\sum_{n=1}^\infty \frac{\cos 2nx}{n^3}$$

$$= \frac{3}{4}x^2 Cl_2(2x) - \frac{3}{8}Cl_4(2x) - \frac{1}{2}x^3 Cl_1(2x) + \frac{3}{4}x Cl_3(2x)$$

Therefore we may write (2.57) as

(2.58) $$\sum_{n=1}^\infty \frac{1}{n^3}\binom{2n}{n}^{-1} \sum_{j=1}^n \frac{(-1)^j}{j}\binom{2n}{n-j}\sin 2jx$$

$$= \frac{1}{3}x^3 \log \sin x + \frac{4}{3}x^3 \log 2 + 2Cl_4(2x) - 2x Cl_3(2x) - 3x\varsigma(3)$$

$$- \frac{1}{4}x^2 Cl_2(2x) + \frac{1}{8}Cl_4(2x) + \frac{1}{3}x^3 Cl_1(2x) - \frac{1}{4}x Cl_3(2x)$$

$$= \frac{1}{3}x^3 \log \sin x + \frac{4}{3}x^3 \log 2 + \frac{17}{8}Cl_4(2x) - \frac{9}{4}x Cl_3(2x) - \frac{1}{4}x^2 Cl_2(2x) + \frac{1}{3}x^3 Cl_1(2x) - 3x\varsigma(3)$$

□

Differentiating (2.1) gives us

(2.59) $$x = \frac{1}{2}\sum_{n=1}^\infty \frac{2^{2n} \sin^{2n-1} x \cos x}{n}\binom{2n}{n}^{-1}$$

and with $x \to \pi x$ this becomes

(2.60) $$\pi x = \frac{1}{2} \sum_{n=1}^{\infty} \frac{2^{2n} \sin^{2n-1} \pi x \cos \pi x}{n} \binom{2n}{n}^{-1}$$

We now multiply this across by $\cot \pi x$ and integrate to obtain

$$\int_0^t \pi x \cot \pi x \, dx = \frac{1}{2} \sum_{n=1}^{\infty} \frac{2^{2n}}{n} \binom{2n}{n}^{-1} \int_0^t \sin^{2n-2} \pi x \cos^2 \pi x \, dx$$

$$= \frac{1}{2} \sum_{n=1}^{\infty} \frac{2^{2n}}{n} \binom{2n}{n}^{-1} \int_0^t \left( \sin^{2n-2} \pi x - \sin^{2n} \pi x \right) dx$$

The Wolfram Mathematica Online Integrator evaluates the above integral in terms of the hypergeometric function $_3F_2\left(\frac{3}{2},\frac{3}{2},-n;\frac{5}{2};\cos^2 \pi t\right)$.

Alternatively, using integration by parts, we get

$$\int_0^t \sin^{2n-2} \pi x \cos^2 \pi x \, dx = \int_0^t \sin^{2n-2} \pi x \cos \pi x \cos \pi x \, dx$$

$$= \cos \pi x \frac{\sin^{2n-1} \pi x}{(2n-1)\pi} \bigg|_0^t + \frac{1}{(2n-1)} \int_0^t \sin^{2n} \pi x \, dx$$

and with $t = 1/2$ we have

$$\int_0^{1/2} \sin^{2n-2} \pi x \cos^2 \pi x \, dx = \frac{1}{2n-1} \int_0^{1/2} \sin^{2n} \pi x \, dx$$

Since

$$\int_0^{1/2} \sin^{2n} \pi x \, dx = \pi \int_0^{\pi/2} \sin^{2n} x \, dx$$

we have

$$\int_0^{\pi/2} \sin^{2n} x \, dx = \frac{(2n-1)!!}{(2n)!!} \frac{\pi}{2} = \frac{1}{2^{2n}} \binom{2n}{n} \frac{\pi}{2}$$

Therefore we see that

$$\int_0^{\frac{1}{2}} \sin^{2n-2}\pi x \cos^2 \pi x \, dx = \frac{1}{2n-1} \frac{1}{2^{2n}} \binom{2n}{n} \frac{\pi^2}{2}$$

Hence we obtain

$$\int_0^{\frac{1}{2}} \pi x \cot \pi x \, dx = \frac{\pi^2}{4} \sum_{n=1}^{\infty} \frac{1}{n(2n-1)}$$

We see that

$$\sum_{n=1}^{\infty} \frac{1}{n(2n-1)} = -\sum_{n=1}^{\infty} \left( \frac{1}{n-\frac{1}{2}} - \frac{1}{n} \right)$$

Since the digamma function may be represented by [43, p.14]

$$\psi(a) = -\gamma - \frac{1}{a} - \sum_{n=1}^{\infty} \left( \frac{1}{n+a} - \frac{1}{n} \right)$$

we see that

$$\sum_{n=1}^{\infty} \left( \frac{1}{n-\frac{1}{2}} - \frac{1}{n} \right) = -\psi\left(-\frac{1}{2}\right) - \gamma + 2$$

We know that [43, p.22]

$$\psi\left(-\frac{1}{2}\right) = 2 - \gamma - 2\log 2$$

and we therefore deduce that

$$\sum_{n=1}^{\infty} \left( \frac{1}{n-\frac{1}{2}} - \frac{1}{n} \right) = 2\log 2$$

Hence we have the well-known integral

$$\int_0^{\frac{1}{2}} x \cot \pi x \, dx = \frac{\pi}{2} \log 2$$

□

We now divide (2.59) by $\sin x$ and integrate to obtain

$$\int_0^t \frac{x}{\sin x}\,dx = \frac{1}{2}\sum_{n=1}^{\infty} \frac{2^{2n}\sin^{2n-1}t}{n(2n-1)}\binom{2n}{n}^{-1}$$

Alternatively, with $n \to n+1$ in the summation this may be written as

(2.61) $$\int_0^t \frac{x}{\sin x}\,dx = \sum_{n=0}^{\infty} \frac{2^{2n}\sin^{2n+1}t}{(2n+1)^2}\binom{2n}{n}^{-1}$$

We recall from [22] that

(2.62) $$\int_a^b \frac{p(x)}{\sin x}\,dx = 2\sum_{n=0}^{\infty}\int_a^b p(x)\sin(2n+1)x\,dx$$

where $p(x)$ is a twice continuously differentiable function. It should be noted that in the above formula we require either (i) both $\sin(x/2)$ and $\cos(x/2)$ have no zero in $[a,b]$ or (ii) if either $\sin(a/2)$ or $\cos(a/2)$ is equal to zero then $p(a)$ must also be zero. Condition (i) is equivalent to the requirement that $\sin x$ has no zero in $[a,b]$.

Letting $p(x) = x$ in (2.62) we get

(2.63) $$\int_0^t \frac{x}{\sin x}\,dx = -2t\sum_{n=0}^{\infty}\frac{\cos(2n+1)t}{2n+1} + 2\sum_{n=0}^{\infty}\frac{\sin(2n+1)t}{(2n+1)^2}$$

and hence we obtain

(2.64) $$\sum_{n=0}^{\infty}\frac{2^{2n}\sin^{2n+1}t}{(2n+1)^2}\binom{2n}{n}^{-1} = -t\sum_{n=0}^{\infty}\frac{\cos(2n+1)t}{2n+1} + \sum_{n=0}^{\infty}\frac{\sin(2n+1)t}{(2n+1)^2}$$

In particular we have from (2.63)

(2.65) $$\int_0^{\pi/2} \frac{x}{\sin x}\,dx = 2\sum_{n=0}^{\infty}\frac{(-1)^n}{(2n+1)^2} = 2G$$

and we then determine that

(2.66) $$G = \frac{1}{2}\sum_{n=0}^{\infty}\frac{2^{2n}}{(2n+1)^2}\binom{2n}{n}^{-1}$$

in agreement with Sherman's compendium of formulae [42, Eq.(3.114)].

Using integration by parts Bradley [13a] showed that

$$\int_0^z \log \tan u \, du - z \log \tan z = -\int_0^z \frac{u \sec^2 u}{\tan u} du$$

$$= -\frac{1}{4}\int_0^{2z} \frac{u \, du}{\tan(u/2)\cos^2(u/2)}$$

$$= -\frac{1}{2}\int_0^{2z} \frac{u}{\sin u} du$$

$$= -\int_0^{\sin(2z)} \frac{2x \sin^{-1} x}{\sqrt{1-x^2}} \frac{dx}{x^2}$$

Differentiating (2.0) gives us

$$\frac{2x \sin^{-1} x}{\sqrt{1-x^2}} = \sum_{n=1}^{\infty} \frac{(2x)^{2n}}{n} \binom{2n}{n}^{-1} \quad \text{for } |x| \leq 1$$

and integrating term by term results in for $0 \leq z \leq \pi/4$

$$\int_0^z \log \tan u \, du - z \log \tan z = -\frac{1}{4}\sum_{n=0}^{\infty} \frac{(2\sin 2z)^{2n+1}}{(2n+1)^2} \binom{2n}{n}^{-1}$$

Therefore we obtain (2.61) again

$$\int_0^z \frac{u}{\sin u} du = \sum_{n=0}^{\infty} \frac{2^{2n} \sin^{2n+1} z}{(2n+1)^2} \binom{2n}{n}^{-1}$$

We note from (2.68) below that

$$\int_0^t \frac{x}{\sin x} dx = t \log \tan\left(\frac{t}{2}\right) + 2\sum_{n=0}^{\infty} \frac{\sin(2n+1)t}{(2n+1)^2}$$

and hence we have

(2.67) $$t \log \tan t + \sum_{n=0}^{\infty} \frac{\sin 2(2n+1)t}{(2n+1)^2} = \sum_{n=0}^{\infty} \frac{2^{2n} \sin^{2n+1} t}{(2n+1)^2} \binom{2n}{n}^{-1}$$

□

Reference to [44, p.149] shows that for $0 < t < \pi$ (this may also be easily derived using the methods outlined in [22])

$$\sum_{n=0}^{\infty} \frac{\cos(2n+1)t}{2n+1} = -\frac{1}{2} \log \tan\left(\frac{t}{2}\right)$$

and for $0 \leq t \leq \pi$

$$\sum_{n=0}^{\infty} \frac{\sin(2n+1)t}{(2n+1)^2} = -\frac{1}{2} \int_0^t \log \tan\left(\frac{x}{2}\right)$$

Referring to (2.63) we have

(2.68) $$\int_0^t \frac{x}{\sin x} dx = t \log \tan\left(\frac{t}{2}\right) + 2 \sum_{n=0}^{\infty} \frac{\sin(2n+1)t}{(2n+1)^2}$$

and noting the Clausen function

$$\text{Cl}_2(t) = -\int_0^t \log[2\sin(u/2)] \, du = \sum_{n=1}^{\infty} \frac{\sin nt}{n^2}$$

we see that

$$\text{Cl}_2(\pi - t) = \sum_{n=0}^{\infty} \frac{\sin n(\pi - t)}{n^2} = \sum_{n=0}^{\infty} (-1)^n \frac{\sin nt}{n^2}$$

$$\text{Cl}_2(t) + \text{Cl}_2(\pi - t) = 2 \sum_{n=0}^{\infty} \frac{\sin(2n+1)t}{(2n+1)^2}$$

Hence, as shown by Lewin [39, p.255] we have

(2.69) $$\int_0^t \frac{x}{\sin x} dx = t \log \tan\left(\frac{t}{2}\right) + \text{Cl}_2(t) + \text{Cl}_2(\pi - t)$$

More generally we have

(2.70) $$\sum_{n=0}^{\infty}\frac{\sin(2n+1)t}{(2n+1)^2} - t\sum_{n=0}^{\infty}\frac{\cos(2n+1)t}{2n+1} = \sum_{n=0}^{\infty}\frac{2^{2n}\sin^{2n+1}t}{(2n+1)^2}\binom{2n}{n}^{-1}$$

$$= t\log\tan\left(\frac{t}{2}\right) + Cl_2(t) + Cl_2(\pi - t)$$

See also the recent paper by Cho et al. [16] which considers the integrals $\int_0^t \frac{x^m}{\sin x}\,dx$ and in particular

$$\int_0^{\pi/3}\frac{x}{\sin x}\,dx = -\frac{\pi}{6}\log 3 - \frac{\sqrt{3}}{9}\pi^2 + \frac{\sqrt{3}}{18}\varsigma\left(2,\frac{1}{6}\right)$$

Such integrals may be easily evaluated using the method set out in [22] (and we have the added advantage that the upper end of the interval of integration does not need to be restricted to the form $\pi/p$ where $p$ is an integer).

$\square$

We have

$$\int_0^t \sin^{2n+1}x\,dx = \int_0^t \sin^{2n}x \sin x\,dx = \int_0^t (1-\cos^2 x)^n \sin x\,dx$$

$$= \int_0^t \sum_{k=0}^n \binom{n}{k}(-1)^k \cos^{2k}x \sin x\,dx$$

and we therefore obtain

$$\int_0^t \sin^{2n+1}x\,dx = \sum_{k=0}^n \binom{n}{k}\frac{(-1)^k}{2k+1} - \sum_{k=0}^n \binom{n}{k}\frac{(-1)^k \cos^{2k+1}t}{2k+1}$$

With $t = \pi/2$ we get

$$\int_0^{\pi/2} \sin^{2n+1}x\,dx = \sum_{k=0}^n \binom{n}{k}\frac{(-1)^k}{2k+1}$$

as compared with the frequently quoted form of the Wallis integral formula

$$\int_0^{\pi/2} \sin^{2n+1}x\,dx = \frac{(2n)!!}{(2n+1)!!}$$

We therefore have the following identity which is reported in [33, p.270]

$$(2.71) \quad \sum_{k=0}^{n}\binom{n}{k}(-1)^k \frac{1}{2k+1} = \frac{(2n)!!}{(2n+1)!!} = \frac{[2^n n!]^2}{(2n+1)!}$$

so that

$$(2.72) \quad \int_0^{\pi/2} \sin^{2n+1} x\, dx = \frac{[2^n n!]^2}{(2n+1)!}$$

We recall (1.5)

$$\sum_{k=0}^{n}\binom{n}{k}\frac{x^k}{(k+y)^s} = \frac{1}{\Gamma(s)}\int_0^\infty u^{s-1} e^{-yu}(1+xe^{-u})^n\, du$$

and note that

$$\sum_{k=0}^{n}\binom{n}{k}(-1)^k \frac{1}{2k+1} = \frac{1}{2}\sum_{k=0}^{n}\binom{n}{k}(-1)^k \frac{1}{k+\frac{1}{2}} = \frac{1}{2}\int_0^\infty e^{-u/2}(1-xe^{-u})^n\, du$$

□

Integrating (2.70) gives us

$$2\sum_{n=0}^{\infty}\frac{1}{(2n+1)^3} - 2\sum_{n=0}^{\infty}\frac{\cos(2n+1)t}{(2n+1)^3} - t\sum_{n=0}^{\infty}\frac{\sin(2n+1)t}{(2n+1)^2} = \sum_{n=0}^{\infty}\frac{2^{2n}}{(2n+1)^2}\binom{2n}{n}^{-1}\int_0^t \sin^{2n+1} x\, dx$$

and with $t = \pi/2$ we have

$$\frac{7}{4}\varsigma(3) + \frac{\pi}{2}\sum_{n=0}^{\infty}\frac{(-1)^{n+1}}{(2n+1)^2} = \sum_{n=0}^{\infty}\frac{2^{2n}}{(2n+1)^2}\binom{2n}{n}^{-1}\int_0^{\pi/2} \sin^{2n+1} x\, dx$$

$$= \sum_{n=0}^{\infty}\frac{2^{2n}}{(2n+1)^2}\binom{2n}{n}^{-1}\frac{[2^n n!]^2}{(2n+1)!}$$

Hence we obtain

$$(2.73) \quad \frac{7}{4}\varsigma(3) - \frac{1}{2}\pi G = \sum_{n=0}^{\infty}\frac{2^{4n}}{(2n+1)^3}\frac{[n!]^4}{[(2n)!]^2}$$

which was also previously obtained by Batir [9].

□

We now multiply (2.1) by $\sin x$ and integrate to obtain

$$\int_0^{\pi/2} x^2 \sin x \, dx = \frac{1}{2} \sum_{n=1}^{\infty} \frac{[n!]^2}{n^2 (2n)!} 2^{2n} \int_0^{\pi/2} \sin^{2n+1} x \, dx$$

and therefore we see that

(2.74) $$\pi - 2 = \frac{1}{2} \sum_{n=1}^{\infty} \frac{2^{4n} [n!]^4}{n^2 (2n+1)[(2n)!]^2}$$

Let us now divide (2.1) by $\sin x$ and integrate to obtain

$$\int_0^{\pi/2} \frac{x^2}{\sin x} dx = \frac{1}{2} \sum_{n=1}^{\infty} \frac{[n!]^2}{n^2 (2n)!} 2^{2n} \int_0^{\pi/2} \sin^{2n-1} x \, dx$$

In equation (6.29) in [22] we have previously shown that

$$\int_0^{\pi/2} \frac{x^2}{\sin x} dx = 2\pi G - \frac{7}{2} \varsigma(3)$$

(the evaluation of this integral is also contained in Bradley's website [13]) and using the integral in [12, p.113]

Using

$$\int_0^{\pi/2} \sin^{2n-1} x \, dx = \sum_{k=0}^{n-1} \binom{n-1}{k} \frac{(-1)^k}{2k+1} = \frac{(2n-2)!!}{(2n-1)!!}$$

we obtain

(2.75) $$8\pi G - 14\varsigma(3) = \sum_{n=1}^{\infty} \frac{2^{4n} [n!]^4}{n^3 [(2n)!]^2}$$

which was previously derived by Batir [9].

□

Letting $p(x) = x \cos x$ in (2.62) gives us

$$\int_0^t \frac{x\cos x}{\sin x}\,dx = \int_0^t x\cot x\,dx = 2\sum_{n=0}^{\infty}\int_0^t x\cos x\sin(2n+1)x\,dx$$

and we have using integration by parts for $n \geq 1$

$$8\int_0^t x\cos x\sin(2n+1)x\,dx = \frac{\sin(2nt)}{n^2} - 2t\frac{\cos(2nt)}{n} + \frac{\sin[2(n+1)t]}{(n+1)^2} - 2t\frac{\cos[2(n+1)t]}{n+1}$$

This gives us

$$2\int_0^{\pi/2} x\cos x\sin(2n+1)x\,dx = \frac{\pi}{4}\left[\frac{(-1)^n}{n+1} - \frac{(-1)^n}{n}\right]$$

and thus

$$2\sum_{n=0}^{\infty}\int_0^{\pi/2} x\cos x\sin(2n+1)x\,dx = \frac{\pi}{2}\log 2$$

This of course concurs with the well-known Euler integral

$$\int_0^{\pi/2} x\cot x\,dx = \frac{\pi}{2}\log 2$$

$\square$

Dividing (2.1) by $\sin^2 x$ and integrating gives us

$$\int_0^t \frac{x^2}{\sin^2 x}\,dx = \frac{1}{2}\sum_{n=1}^{\infty}\frac{[n!]^2}{n^2(2n)!}2^{2n}\int_0^t \sin^{2n-2} x\,dx$$

and with $t = \pi/2$ we have after dealing separately with the case for $n = 1$

$$\int_0^{\pi/2} \frac{x^2}{\sin^2 x}\,dx = \frac{\pi}{2} + \frac{1}{2}\sum_{n=2}^{\infty}\frac{[n!]^2}{n^2(2n)!}2^{2n}\frac{\pi}{2^{2n-1}}\binom{2n-2}{n-1}$$

$$= \frac{\pi}{2} + \pi\sum_{n=2}^{\infty}\frac{1}{(2n-1)(2n-2)}$$

$$= \frac{\pi}{2} + \pi\sum_{n=2}^{\infty}\left(\frac{1}{2n-1} - \frac{1}{2n-2}\right)$$

$$= \frac{\pi}{2} + \pi \sum_{n=1}^{\infty} \left( \frac{1}{2n+1} - \frac{1}{2n} \right)$$

$$= \frac{\pi}{2} + \frac{\pi}{2} \sum_{n=1}^{\infty} \left( \frac{1}{n+\frac{1}{2}} - \frac{1}{n} \right)$$

Since the digamma function may be represented by [43, p.14]

$$\psi(a) = -\gamma - \frac{1}{a} - \sum_{n=1}^{\infty} \left( \frac{1}{n+a} - \frac{1}{n} \right)$$

we see that

$$\sum_{n=1}^{\infty} \left( \frac{1}{n+\frac{1}{2}} - \frac{1}{n} \right) = -\psi\left(\frac{1}{2}\right) - \gamma - 2$$

We know that [43, p.20]

$$\psi\left(\frac{1}{2}\right) = -\gamma - 2\log 2$$

and we therefore deduce that

$$\int_0^{\pi/2} \frac{x^2}{\sin^2 x} dx = \pi \log 2$$

Using integration by parts we see that

$$\int_0^{\pi/2} \frac{x^2}{\sin^2 x} dx = 2 \int_0^{\pi/2} x \cot x \, dx$$

and thus we have

$$\int_0^{\pi/2} x \cot x \, dx = -\int_0^{\pi/2} \log \sin x \, dx$$

Euler [25] was the first person to show that

$$\int_0^{\pi/2} \log \sin x \, dx = -\frac{\pi}{2} \log 2$$

and we have therefore come full circle! □

It may also be noted that Ramanujan [11, Part I, p.261] showed that for $|x| < 2\pi$

$$(2.76) \quad \frac{1}{2}\int_0^x u^n \cot\left(\frac{u}{2}\right) du = \cos\left(\frac{n\pi}{2}\right) n!\varsigma(n+1) - \sum_{j=0}^n (-1)^{j(j+1)/2} \frac{\Gamma(n+1)}{\Gamma(n+1-j)} x^{n-j} \text{Cl}_{j+1}(x)$$

and also see the paper by Srivastava et al [43a, Eq (5.12)]. Ramanujan's formula was recently employed by Muzaffar [40]. With $x = \pi$ we have

$$(2.77) \quad 2^n \int_0^{\pi/2} x^n \cot x \, dx = \cos\left(\frac{n\pi}{2}\right) n!\varsigma(n+1) - \sum_{j=0}^n (-1)^{j(j+1)/2} \frac{\Gamma(n+1)}{\Gamma(n+1-j)} \pi^{n-j} \text{Cl}_{j+1}(\pi)$$

□

Provided $a \neq b$ we readily determine that

$$\int \sin ax \sin bx \, dx = \frac{1}{2}\left[\frac{\sin(a-b)x}{a-b} - \frac{\sin(a+b)x}{a+b}\right] + c$$

and hence we have

$$\int_0^t \sin px \sin(2n+1)x \, dx = \frac{1}{2}\left[\frac{\sin(p-2n-1)t}{p-2n-1} - \frac{\sin(p+2n+1)t}{p+2n+1}\right]$$

and more specifically

$$\int_0^{\pi/2} \sin px \sin(2n+1)x \, dx = \frac{(-1)^{n+1} \cos(p\pi/2)}{2}\left[\frac{1}{p-2n-1} + \frac{1}{p+2n+1}\right]$$

Therefore with $p(x) = \sin px$ in (2.62) we obtain

$$\int_0^{\pi/2} \frac{\sin px}{\sin x} dx = 2\sum_{n=0}^{\infty} \int_0^{\pi/2} \sin px \sin(2n+1)x \, dx$$

so that

(2.78) $$\int_0^{\pi/2} \frac{\sin px}{\sin x} dx = 2p \cos(p\pi/2) \sum_{n=0}^{\infty} \frac{(-1)^n}{(2n+1)^2 - p^2}$$

and, using L'Hôpital's rule, we see that in the limit as $p \to 1$

$$\int_0^{\pi/2} dx = \lim_{p \to 1} \frac{-p\pi \sin(p\pi/2) + 2\cos(p\pi/2)}{-2p} = \frac{\pi}{2}$$

With $p = 2$ we obtain

(2.79) $$\sum_{n=0}^{\infty} \frac{(-1)^{n+1}}{(2n+1)^2 - 4^2} = \frac{1}{2}$$

We may write (2.78) as

$$\sum_{n=0}^{\infty} \frac{(-1)^n}{(2n+1)^2 - p^2} = \frac{1}{2p \cos(p\pi/2)} \int_0^{\pi/2} \frac{\sin px}{\sin x} dx$$

and, using L'Hôpital's rule, we obtain as $p \to 0$

$$\sum_{n=0}^{\infty} \frac{(-1)^n}{(2n+1)^2} = \lim_{p \to 0} \left[ \frac{1}{-p\pi \sin(p\pi/2) + 2\cos(p\pi/2)} \int_0^{\pi/2} \frac{x \cos px}{\sin x} dx \right]$$

Hence we obtain (2.65) again

$$2 \sum_{n=0}^{\infty} \frac{(-1)^n}{(2n+1)^2} = \int_0^{\pi/2} \frac{x}{\sin x} dx$$

We note from Fettis [26b] that

(2.80) $$\int_0^{\pi/2} \frac{\sin px}{\sin x} dx = \frac{\pi}{2} - \frac{1}{2} \cos\left(\frac{p\pi}{2}\right) \left[ \psi\left(\frac{3+p}{4}\right) - \psi\left(\frac{1+p}{4}\right) \right]$$

and

(2.81) $$\int_0^{\pi/2} \frac{1 - \cos px}{\sin x} dx = \psi\left(\frac{1+p}{2}\right) - \psi\left(\frac{1}{2}\right) - \frac{1}{2} \sin\left(\frac{p\pi}{2}\right) \left[ \psi\left(\frac{3+p}{4}\right) - \psi\left(\frac{1+p}{4}\right) \right]$$

Comparing (2.78) with (2.80) gives us

(2.82) $$2p\cos(p\pi/2)\sum_{n=0}^{\infty}\frac{(-1)^n}{(2n+1)^2-p^2}=\frac{\pi}{2}-\frac{1}{2}\cos\left(\frac{p\pi}{2}\right)\left[\psi\left(\frac{3+p}{4}\right)-\psi\left(\frac{1+p}{4}\right)\right]$$

Differentiating (2.80) results in

$$\int_0^{\pi/2}\frac{x\cos px}{\sin x}dx=-\frac{1}{8}\cos\left(\frac{p\pi}{2}\right)\left[\psi'\left(\frac{3+p}{4}\right)-\psi'\left(\frac{1+p}{4}\right)\right]+\frac{\pi}{4}\sin\left(\frac{p\pi}{2}\right)\left[\psi\left(\frac{3+p}{4}\right)-\psi\left(\frac{1+p}{4}\right)\right]$$

and with $p=0$ we obtain

$$\int_0^{\pi/2}\frac{x}{\sin x}dx=-\frac{1}{8}\left[\psi'\left(\frac{3}{4}\right)-\psi'\left(\frac{1}{4}\right)\right]$$

Using (2.65) we see that

(2.83) $$2G=-\frac{1}{8}\left[\psi'\left(\frac{3}{4}\right)-\psi'\left(\frac{1}{4}\right)\right]$$

which concurs with Kölbig [33a].

With $p=1$ we have the well known integral

$$\int_0^{\pi/2}x\cot x\,dx=\frac{\pi}{4}\left[\psi(1)-\psi\left(\frac{1}{2}\right)\right]=\frac{\pi}{2}\log 2$$

With $p=2$ we obtain

$$\int_0^{\pi/2}\left[\frac{x}{\sin x}-2x\sin x\right]dx=\frac{1}{8}\left[\psi'\left(1+\frac{1}{4}\right)-\psi'\left(\frac{3}{4}\right)\right]$$

and since $\psi'\left(1+\frac{1}{4}\right)=\psi'\left(\frac{1}{4}\right)-16$ we obtain (2.83) again.

Differentiating (2.81) results in

$$\int_0^{\pi/2}\frac{x\sin px}{\sin x}dx=\frac{1}{2}\psi'\left(\frac{1+p}{2}\right)-\frac{1}{8}\sin\left(\frac{p\pi}{2}\right)\left[\psi'\left(\frac{3+p}{4}\right)-\psi'\left(\frac{1+p}{4}\right)\right]$$

$$-\frac{\pi}{4}\cos\left(\frac{p\pi}{2}\right)\left[\psi\left(\frac{3+p}{4}\right)-\psi\left(\frac{1+p}{4}\right)\right]$$

and with $p=0$ we obtain

$$\frac{1}{2}\psi'\left(\frac{1}{2}\right)=\frac{\pi}{4}\left[\psi\left(\frac{3}{4}\right)-\psi\left(\frac{1}{4}\right)\right]$$

so that

(2.84) $\quad \psi'\left(\frac{1}{2}\right)=\frac{\pi^2}{2}$

in accordance with Kölbig's paper [33a].

With $p=1$ we obtain

$$\frac{\pi^2}{8}=\frac{1}{2}\psi'(1)-\frac{1}{8}\left[\psi'(1)-\psi'\left(\frac{1}{2}\right)\right]$$

and, since $\psi'(1)=\varsigma(2)$, this concurs with (2.84).

We have provided $p \neq 2n+1$

$$\int_0^t \cos px \sin(2n+1)x\, dx = \frac{1}{2}\left[\frac{\cos(p-2n-1)t}{p-2n-1}-\frac{\cos(p+2n+1)t}{p+2n+1}\right]-\frac{2n+1}{p^2-(2n+1)^2}$$

and thus

$$\int_0^{\pi/2} \cos px \sin(2n+1)x\, dx = \frac{(-1)^{n+1}}{2}\left[\frac{1}{p-2n-1}+\frac{1}{p+2n+1}\right]-\frac{2n+1}{p^2-(2n+1)^2}$$

$$=\frac{(-1)^{n+1}p-(2n+1)}{p^2-(2n+1)^2}$$

Therefore letting $p(x)=1-\cos px$ in (2.62) gives us

(2.85) $\quad \displaystyle\int_0^{\pi/2}\frac{1-\cos px}{\sin x}dx = 2\sum_{n=0}^{\infty}\left[\frac{1}{2n+1}+\frac{(-1)^n p+(2n+1)}{p^2-(2n+1)^2}\right]$

and comparing this with (2.81) results in

$$2\sum_{n=0}^{\infty}\left[\frac{1}{2n+1}+\frac{(-1)^n p+(2n+1)}{p^2-(2n+1)^2}\right]=\psi\left(\frac{1+p}{2}\right)-\psi\left(\frac{1}{2}\right)-\frac{1}{2}\sin\left(\frac{p\pi}{2}\right)\left[\psi\left(\frac{3+p}{4}\right)-\psi\left(\frac{1+p}{4}\right)\right]$$

Multiplying across by $\cos\left(\frac{p\pi}{2}\right)$ gives us

$$2\cos\left(\frac{p\pi}{2}\right)\sum_{n=0}^{\infty}\left[\frac{1}{2n+1}+\frac{(-1)^n p+(2n+1)}{p^2-(2n+1)^2}\right]$$

$$=\cos\left(\frac{p\pi}{2}\right)\left[\psi\left(\frac{1+p}{2}\right)-\psi\left(\frac{1}{2}\right)\right]-\frac{1}{2}\cos\left(\frac{p\pi}{2}\right)\sin\left(\frac{p\pi}{2}\right)\left[\psi\left(\frac{3+p}{4}\right)-\psi\left(\frac{1+p}{4}\right)\right]$$

We have

$$2\cos\left(\frac{p\pi}{2}\right)\sum_{n=0}^{\infty}\left[\frac{1}{2n+1}+\frac{(-1)^n p+(2n+1)}{p^2-(2n+1)^2}\right]$$

$$=\frac{\pi}{2}-\frac{1}{2}\cos\left(\frac{p\pi}{2}\right)\left[\psi\left(\frac{3+p}{4}\right)-\psi\left(\frac{1+p}{4}\right)\right]+2\cos\left(\frac{p\pi}{2}\right)\sum_{n=0}^{\infty}\left[\frac{1}{2n+1}+\frac{2n+1}{p^2-(2n+1)^2}\right]$$

$$=\frac{\pi}{2}-\frac{1}{2}\cos\left(\frac{p\pi}{2}\right)\left[\psi\left(\frac{3+p}{4}\right)-\psi\left(\frac{1+p}{4}\right)\right]+2p^2\cos\left(\frac{p\pi}{2}\right)\sum_{n=0}^{\infty}\left[\frac{1}{2n+1}\frac{1}{p^2-(2n+1)^2}\right]$$

and therefore we see that

(2.86) $$2p^2\sum_{n=0}^{\infty}\frac{1}{2n+1}\frac{1}{p^2-(2n+1)^2}$$

$$=\psi\left(\frac{1+p}{2}\right)-\psi\left(\frac{1}{2}\right)+\frac{1}{2}\left[1-\sin\left(\frac{p\pi}{2}\right)\right]\left[\psi\left(\frac{3+p}{4}\right)-\psi\left(\frac{1+p}{4}\right)\right]-\frac{\pi}{2}/\cos\left(\frac{p\pi}{2}\right)$$

This is reminiscent of the expression in Prudnikov et al. [41a, 5.1.25-13]

(2.87) $$2p^2\sum_{n=1}^{\infty}\frac{1}{n}\frac{1}{p^2-n^2}=\psi(1+p)+\psi(1-p)+2\gamma$$

and separating the even and odd terms gives us

$$\sum_{n=1}^{\infty}\frac{1}{n}\frac{1}{p^2-n^2} = \sum_{n=1}^{\infty}\frac{1}{2n}\frac{1}{p^2-4n^2} + \sum_{n=0}^{\infty}\frac{1}{2n+1}\frac{1}{p^2-(2n+1)^2}$$

Letting $p \to p/2$ in (2.87) we see that

$$\psi\left(1+\frac{p}{2}\right)+\psi\left(1-\frac{p}{2}\right)+2\gamma = 2p^2\sum_{n=1}^{\infty}\frac{1}{n}\frac{1}{p^2-4n^2}$$

and hence we obtain

$$\psi(1+p)+\psi(1-p)+2\gamma - \frac{1}{2}\left[\psi\left(1+\frac{p}{2}\right)+\psi\left(1-\frac{p}{2}\right)+2\gamma\right]$$

$$= \psi\left(\frac{1+p}{2}\right)-\psi\left(\frac{1}{2}\right)+\frac{1}{2}\left[1-\sin\left(\frac{p\pi}{2}\right)\right]\left[\psi\left(\frac{3+p}{4}\right)-\psi\left(\frac{1+p}{4}\right)\right]-\frac{\pi}{2}/\cos\left(\frac{p\pi}{2}\right)$$

which may be written as

(2.88) $\psi(1+p)+\psi(1-p)-\frac{1}{2}\left[\psi\left(1+\frac{p}{2}\right)+\psi\left(1-\frac{p}{2}\right)\right]$

$$= \psi\left(\frac{1+p}{2}\right)+2\log 2 +\frac{1}{2}\left[1-\sin\left(\frac{p\pi}{2}\right)\right]\left[\psi\left(\frac{3+p}{4}\right)-\psi\left(\frac{1+p}{4}\right)\right]-\frac{\pi}{2}/\cos\left(\frac{p\pi}{2}\right)$$

For example, with $p=1/2$ we obtain

$$\psi\left(1+\frac{1}{2}\right)+\psi\left(\frac{1}{2}\right)-\frac{1}{2}\left[\psi\left(1+\frac{1}{4}\right)+\psi\left(\frac{3}{4}\right)\right]$$

$$= \psi\left(\frac{3}{4}\right)+2\log 2 + \frac{1}{2}\left[1-\frac{1}{2}\sqrt{2}\right]\left[\psi\left(\frac{7}{8}\right)-\psi\left(\frac{3}{8}\right)\right]-\frac{\pi\sqrt{2}}{2}$$

and this may be verified by substituting the specific values for the digamma function in [43, p.20] (albeit one of the relevant values is reported incorrectly in [43, p.20]).

□

Ramanujan [11, Part I, p.263] reported the following Maclaurin series (in the slightly modified form employed by Borwein and Chamberland [12a])

(2.89) $$\left(\sin^{-1} y\right)^4 = \frac{3}{2}\sum_{n=1}^{\infty} \frac{H_{n-1}^{(2)}}{n^2}\binom{2n}{n}^{-1}(2y)^{2n} \qquad , -1 \leq y \leq 1$$

where $H_n^{(m)}$ is the generalised harmonic number

$$H_n^{(m)} = \sum_{k=1}^{n} \frac{1}{k^m}$$

Therefore upon letting $x = \sin^{-1} y$ we get

(2.90) $$x^4 = \frac{3}{2}\sum_{n=1}^{\infty} \frac{H_{n-1}^{(2)}}{n^2}\binom{2n}{n}^{-1} 2^{2n} \sin^{2n} x \qquad , -\pi/2 \leq x \leq \pi/2$$

Integration results in

$$\frac{\pi^5}{160} = \frac{3}{2}\sum_{n=1}^{\infty} \frac{H_{n-1}^{(2)}}{n^2}\binom{2n}{n}^{-1} 2^{2n} \int_0^{\pi/2} \sin^{2n} x\, dx$$

and hence we obtain the Euler sum

$$\frac{\pi^5}{160} = \frac{3\pi}{4}\sum_{n=1}^{\infty} \frac{H_{n-1}^{(2)}}{n^2}$$

or equivalently

$$\frac{\pi^4}{120} = \sum_{n=1}^{\infty} \frac{H_{n-1}^{(2)}}{n^2} = \sum_{n=1}^{\infty} \frac{H_n^{(2)}}{n^2} - \varsigma(4)$$

Since $\varsigma(4) = \frac{\pi^4}{90}$ this may be written in its more familiar form as (see for example [19])

$$\frac{7}{4}\varsigma(4) = \sum_{n=1}^{\infty} \frac{H_n^{(2)}}{n^2}$$

Alternatively we now multiply equation (2.90) by $\cot x$ and integrate to obtain

$$\int_0^t x^4 \cot x\, dx = \frac{3}{2}\sum_{n=1}^{\infty} \frac{H_{n-1}^{(2)}}{n^2}\binom{2n}{n}^{-1} 2^{2n} \int_0^t \sin^{2n-1} x \cos x\, dx$$

which results in

(2.91) $$\int_0^t x^4 \cot x\, dx = \frac{3}{4}\sum_{n=1}^{\infty} \frac{2^{2n} H_{n-1}^{(2)} \sin^{2n} t}{n^3} \binom{2n}{n}^{-1}$$

This integral may also be evaluated in terms of, for example, the Clausen functions by using (2.4).

We multiply (2.90) by $\cos x$ and integrate to obtain for $-\pi/2 \leq x \leq \pi/2$

$$4t(t^2 - 6)\cos t + (t^4 - 12t^2 + 24)\sin t = \frac{3}{2}\sum_{n=1}^{\infty} \frac{2^{2n} H_{n-1}^{(2)}}{(2n+1)n^2}\binom{2n}{n}^{-1} \sin^{2n+1} t$$

and with $t = \pi/2$ we obtain

(2.92) $$\frac{\pi^4}{16} - 3\pi^2 + 24 = \frac{3}{2}\sum_{n=1}^{\infty} \frac{2^{2n} H_{n-1}^{(2)}}{(2n+1)n^2}\binom{2n}{n}^{-1}$$

□

We showed in equation (4.3.158) of [19] that for $0 \leq t < 1$

(2.93) $$\int_0^t \pi x \cot \pi x\, dx = \varsigma'(-1, t) - \varsigma'(-1, 1-t) + t\log(2\sin \pi t)$$

and we now refer to the well-known formula [12, p.130]

(2.94) $$\pi x \cot \pi x = -2\sum_{n=0}^{\infty} \varsigma(2n) x^{2n} \qquad , (|x| < 1)$$

(where $\varsigma(0) = -1/2$). This gives us

(2.95) $$-2\sum_{n=0}^{\infty} \frac{\varsigma(2n)}{2n+1} t^{2n+1} = \varsigma'(-1, t) - \varsigma'(-1, 1-t) + t\log(2\sin \pi t)$$

We also have [43, p.12] in terms of the digamma function

$$\pi x \cot \pi x = x\psi(1-x) - x\psi(x)$$

and therefore we have

$$\int_0^t [x\psi(1-x) - x\psi(x)]dx = \varsigma'(-1, t) - \varsigma'(-1, 1-t) + t\log(2\sin \pi t)$$

Integration by parts results in

$$\int_0^t [x\psi(1-x) - x\psi(x)]dx = -t\log[\Gamma(t)\Gamma(1-t)] + \int_0^t \log[\Gamma(x)\Gamma(1-x)]dx$$

whereupon using Euler's reflection formula for the gamma function [43, p.3]

$$\Gamma(t)\Gamma(1-t) = \frac{\pi}{\sin \pi t}$$

this becomes

$$= t\log \sin \pi t - \int_0^t \log \sin \pi x\, dx$$

Therefore, as noted in equation (4.3.158a) in [19], we have

(2.96) $\quad \int_0^t \log[2\sin \pi x]dx = -[\varsigma'(-1,t) - \varsigma'(-1,1-t)]$

which could of course have been obtained more directly by using integration by parts on equation (2.93). This incidentally shows us that

$$\varsigma'(-1,0) = \varsigma'(-1,1) = \varsigma'(-1)$$

Letting $t = 1/2$ in (2.96) immediately gives us Euler's integral

$$\int_0^{1/2} \log \sin \pi x\, dx = -\frac{1}{2}\log 2$$

□

Let us now differentiate (2.95) to obtain

$$-2\sum_{n=0}^{\infty} \varsigma(2n)t^{2n} = \frac{d}{dt}[\varsigma'(-1,t) - \varsigma'(-1,1-t)] + \log(2\sin \pi t) + \pi \cot \pi t$$

Since the Hurwitz zeta function is analytic in the whole complex plane except for $s \neq 1$, its partial derivatives commute in the region where the function is analytic: we therefore have

$$\frac{\partial}{\partial t}\frac{\partial}{\partial s}\varsigma(s,t) = \frac{\partial}{\partial s}\frac{\partial}{\partial t}\varsigma(s,t) = -\frac{\partial}{\partial s}[s\varsigma(s+1,t)]$$

and we then see that

(2.97) $$\frac{\partial}{\partial t}\frac{\partial}{\partial s}\varsigma(s,t) = -\varsigma(s+1,t) - s\frac{\partial}{\partial s}\varsigma(s+1,t)$$

Hence with $s = -1$ we see that

$$\frac{d}{dt}\varsigma'(-1,t) = \varsigma'(0,t) - \varsigma(0,t)$$

$$\frac{d}{dt}\varsigma'(-1,1-t) = -\varsigma'(0,1-t) + \varsigma(0,1-t)$$

and we then obtain

$$-2\sum_{n=0}^{\infty}\varsigma(2n)t^{2n} = [\varsigma'(0,t) + \varsigma'(0,1-t)] - [\varsigma(0,t) - \varsigma(0,1-t)] + \log(2\sin\pi t) + \pi\cot\pi t$$

Using Lerch's identity [10] for $t > 0$

$$\varsigma'(0,t) = \log\Gamma(t) - \frac{1}{2}\log(2\pi)$$

and the well-known result [7, p.264] involving the Bernoulli polynomials $B_n(t)$

$$\varsigma(1-n,t) = -\frac{B_n(t)}{n} \quad \text{for } n \geq 1$$

this becomes

$$-2\sum_{n=0}^{\infty}\varsigma(2n)t^{2n} = \log[\Gamma(t)\Gamma(1-t)] - \log(2\pi) + \log(2\sin\pi t) + \pi\cot\pi t$$

Using Euler's reflection formula for the gamma function

$$\Gamma(t)\Gamma(1-t) = \frac{\pi}{\sin\pi t}$$

we return to where we started from

$$-2\sum_{n=0}^{\infty}\varsigma(2n)t^{2n} = \pi t\cot\pi t$$

Dividing this by $t$, dealing separately with the $n = 0$ term, and integrating results in

$$\sum_{n=1}^{\infty} \varsigma(2n)\frac{t^{2n}}{n} = \log t - \log \sin \pi t + c$$

and in the limit as $t \to 0$ we see that the integration constant is $c = \log \pi$. We thus obtain the known result [43, p.161]

$$\sum_{n=1}^{\infty} \varsigma(2n)\frac{t^{2n}}{n} = \log(\pi t) - \log \sin \pi t$$

We now multiply this by $t$ and integrate to obtain

$$\sum_{n=0}^{\infty} \frac{\varsigma(2n)}{n(2n+2)} u^{2n+2} = u\log(\pi u) - u - \int_0^u t \log \sin \pi t \, dt$$

and, in the case $u = 1/2$, we obtain using (2.4)

$$\varsigma(3) = \frac{2\pi^2}{7}\left[\log \pi - \frac{1}{2} - \sum_{n=1}^{\infty}\frac{\varsigma(2n)}{n(n+1)2^{2n}}\right]$$

This known result was also recently derived by Fujii and Suzuki [28] where, in what appears to be a new approach, they employed the logarithmic form of Euler's infinite product identity for the sine function

$$\log \sin x = \log x + \sum_{n=1}^{\infty} \log\left[1 - \frac{x^2}{n^2\pi^2}\right]$$

□

Upon integrating (2.97) with respect to $t$ we see that

$$-s\int_0^v \varsigma'(s+1,t)dt = \int_0^v \frac{\partial}{\partial t}\frac{\partial}{\partial s}\varsigma(s,t)dt + \int_0^v \varsigma(s+1,t)dt$$

We therefore get

$$-s\int_0^v \varsigma'(s+1,t)dt = \varsigma'(s,v) - \varsigma'(s,0) + \int_0^v \varsigma(s+1,t)dt$$

and with $s = -n$ we have

$$n\int_0^v \varsigma'(1-n,t)dt = \varsigma'(-n,v) - \varsigma'(-n,0) + \int_0^v \varsigma(1-n,t)dt$$

Then using the well-known formula

$$\varsigma(1-n,v) = -\frac{B_n(v)}{n} \text{ for } n \geq 1$$

we obtain

(2.98) $$n\int_0^v \varsigma'(1-n,t)\,dt = \frac{B_{n+1} - B_{n+1}(v)}{n(n+1)} + \varsigma'(-n,v) - \varsigma'(-n,0)$$

This identity was originally derived by Adamchik [1] in a different manner in 1998. With $n = 2$ we obtain

(2.99) $$\int_0^v \varsigma'(-1,t)\,dt = -\frac{1}{12}B_3(v) + \frac{1}{2}\varsigma'(-2,v) + \frac{\varsigma(3)}{8\pi^2}$$

since $\varsigma'(-n,0) = \varsigma'(-n)$ and $\varsigma'(-2) = -\frac{\varsigma(3)}{4\pi^2}$.

With $v = 1$ we note that

(2.100) $$\int_0^1 \varsigma'(1-n,t)\,dt = 0$$

which may also be obtained by integrating (2.29).

## 3. A brief survey of multiple sine functions

The identity (2.6) was also found by Koyama and Kurokawa [34] using the triple sine function. The multiple sine functions are defined for $r = 2, 3, \cdots$ by

(3.1) $$S_r(x) = \exp\left(\frac{x^{r-1}}{r-1}\right) \prod_{n=-\infty, n\neq 0}^{\infty} \left[P_r\left(\frac{x}{n}\right)\right]^{n^{r-1}}$$

$$= \exp\left(\frac{x^{r-1}}{r-1}\right) \prod_{n=1}^{\infty} \left[P_r\left(\frac{x}{n}\right) P_r\left(-\frac{x}{n}\right)\right]^{n^{r-1}}$$

where

$$P_r(u) = (1-u)\exp\left(\frac{u}{1} + \frac{u^2}{2} + \cdots + \frac{u^r}{r}\right)$$

For example, the triple sine function is defined by

(3.2) $$S_3(x) = e^{x^2/2} \prod_{n=1}^{\infty} \left[ \left(1 - \frac{x^2}{n^2}\right)^{n^2} e^{x^2} \right]$$

and we then have

$$\log S_3(x) = \frac{1}{2}x^2 + \sum_{n=1}^{\infty}\left[ n^2 \log\left(1 - \frac{x^2}{n^2}\right) + x^2 \right]$$

whereupon differentiation results in

$$\frac{S_3'(x)}{S_3(x)} = x + \sum_{n=1}^{\infty}\left[ \frac{2n^2 x}{x^2 - n^2} + 2x \right]$$

$$= x + x^2 \sum_{n=1}^{\infty} \frac{2x}{x^2 - n^2}$$

$$= x^2 \left( \frac{1}{x} + \sum_{n=1}^{\infty} \frac{2x}{x^2 - n^2} \right)$$

We have the well-known decomposition formula [12, p.131]

(3.3) $$\pi \cot \pi x = \frac{1}{x} + \sum_{n=1}^{\infty} \frac{2x}{x^2 - n^2}$$

and, since $S_3(0) = 1$, we then see that

(3.4) $$\log S_3(x) = \int_0^x \pi t^2 \cot \pi t \, dt$$

The double sine function is defined by

(3.5) $$S_2(x) = e^x \prod_{n=1}^{\infty}\left[ \left( \frac{1 - x/n}{1 + x/n} \right)^n e^{2x} \right]$$

It is easily seen that

(3.6) $$\log S_2(x) = x + \sum_{n=1}^{\infty}\left[ n\left\{ \log\left(1 - \frac{x}{n}\right) - \log\left(1 - \frac{x}{n}\right) \right\} + 2x \right]$$

and, in the same manner as before, we easily find that

(3.7) $$\log S_2(x) = \int_0^x \pi t \cot \pi t \, dt$$

The function $S_1(x)$ is defined by

(3.8) $$S_1(x) = 2\pi x \prod_{n=1}^{\infty} \left[\left(1 - \frac{x^2}{n^2}\right)\right] = 2\sin \pi x$$

and it is well known that logarithmic differentiation of this results in (3.3) above.

The Barnes double gamma function $\Gamma_2(x) = 1/G(x)$ is defined, inter alia, by [43, p.25]

(3.9) $$G(1+x) = (2\pi)^{x/2} \exp\left[-\frac{1}{2}(\gamma x^2 + x^2 + x)\right] \prod_{k=1}^{\infty} \left\{\left(1 + \frac{x}{k}\right)^k \exp\left(\frac{x^2}{2k} - x\right)\right\}$$

We showed in [23] that the Barnes double gamma function could be represented by

(3.10) $$\log G(1+x) = \frac{1}{2}x\log(2\pi) - \frac{1}{2}x(1+x) - \frac{1}{2}\gamma x^2 = \sum_{n=1}^{\infty}\left[\frac{1}{2n}x^2 - x + n\log\left(1 + \frac{x}{n}\right)\right]$$

and letting $x \to -x$ gives us

$$\log G(1-x) = -\frac{1}{2}x\log(2\pi) + \frac{1}{2}x(1-x) - \frac{1}{2}\gamma x^2 = \sum_{n=1}^{\infty}\left[\frac{1}{2n}x^2 + x + n\log\left(1 - \frac{x}{n}\right)\right]$$

Subtraction results in

$$\log G(1-x) - \log G(1+x)$$

$$= x - x\log(2\pi) + \sum_{n=1}^{\infty}\left[2x + n\log\left(1 - \frac{x}{n}\right) - n\log\left(1 + \frac{x}{n}\right)\right]$$

and, using (3.6), we therefore see that

(3.11) $$\log S_2(x) = \log G(1-x) - \log G(1+x) + x\log(2\pi)$$

We have therefore rediscovered the well-known formula originally found by Kinkelin in 1860 [43, p.30]

(3.12) $$\log \frac{G(1+x)}{G(1-x)} = x\log(2\pi) - \int_0^x \pi t \cot \pi t \, dt$$

which was generalised in 1992 by Freund and Zabrodin [27] who reported the more general identity for $n \geq 2$

$$(3.13) \qquad \Gamma_n(u)[\Gamma_n(-u)]^{(-1)^{n+1}} = \exp\left[-\pi \int_0^u x^{n-1} \cot \pi x \, dx\right]$$

where we have followed the Vignéras notation [43, p.39]

$$\Gamma_n(x) = [G_n(x)]^{(-1)^{n-1}}$$

and $\qquad \Gamma_1(x) = G_1(x) = \Gamma(x) \qquad \Gamma_2(x) = 1/G_2(x) = 1/G(x)$

Koyama and N. Kurokawa [35] have shown that

$$(3.14) \qquad \log S_n(x) = \int_0^x \pi t^{n-1} \cot \pi t \, dt$$

Reference should also be made to Kurokawa's paper [36].

It is an exercise in Bromwich's book [15, p.526] to show that

$$\sum_{n=1}^{\infty}\left[n \log\left(\frac{2n+1}{2n-1}\right) - 1\right] = \frac{1}{2}(1 - \log 2)$$

and it may be noted that this agrees with (3.6) when $x = 1/2$.

$\square$

The gamma function may be defined by [43, p.2]

$$\log \Gamma(x+1) = \sum_{n=1}^{\infty}\left[x \log\left(1 + \frac{1}{n}\right) - \log\left(1 + \frac{x}{n}\right)\right]$$

and this may be written as

$$\log \Gamma(x+1) = \sum_{n=1}^{\infty}\left[x \log\left(1 + \frac{1}{n}\right) + \sum_{k=1}^{\infty}\frac{(-1)^k x^k}{kn^k}\right]$$

$$= \sum_{n=1}^{\infty}\left[x\left[\log\left(1 + \frac{1}{n}\right) - \frac{1}{n}\right] + \sum_{k=2}^{\infty}\frac{(-1)^k x^k}{kn^k}\right]$$

$$= \sum_{n=1}^{\infty}\left[\log\left(1+\frac{1}{n}\right)-\frac{1}{n}\right]x + \sum_{n=1}^{\infty}\sum_{k=2}^{\infty}\frac{(-1)^k x^k}{kn^k}$$

$$= \sum_{n=1}^{\infty}\left[\log\left(1+\frac{1}{n}\right)-\frac{1}{n}\right]x + \sum_{k=2}^{\infty}\frac{(-1)^k x^k}{k}\sum_{n=1}^{\infty}\frac{1}{n^k}$$

Hence we have the well-known Maclaurin series [43, p.160] for $\log\Gamma(x+1)$

(3.15) $$\log\Gamma(x+1) = -\gamma x + \sum_{k=2}^{\infty}\frac{(-1)^k \varsigma(k)}{k}x^k$$

Now referring back to (3.10)

$$\log G(1+x) = \frac{1}{2}x\log(2\pi) - \frac{1}{2}x(1+x) - \frac{1}{2}\gamma x^2 + \sum_{n=1}^{\infty}\left[\frac{1}{2n}x^2 - x + n\log\left(1+\frac{x}{n}\right)\right]$$

and applying the same procedure as with $\log\Gamma(x+1)$ above this becomes

$$= \frac{1}{2}x\log(2\pi) - \frac{1}{2}x(1+x) - \frac{1}{2}\gamma x^2 + \sum_{n=1}^{\infty}\left[\frac{1}{2n}x^2 - x - \sum_{k=1}^{\infty}\frac{(-1)^k x^k}{kn^{k-1}}\right]$$

$$= \frac{1}{2}x\log(2\pi) - \frac{1}{2}x(1+x) - \frac{1}{2}\gamma x^2 - \sum_{n=1}^{\infty}\sum_{k=3}^{\infty}\frac{(-1)^k x^k}{kn^{k-1}}$$

$$= \frac{1}{2}x\log(2\pi) - \frac{1}{2}x(1+x) - \frac{1}{2}\gamma x^2 - \sum_{n=1}^{\infty}\sum_{m=2}^{\infty}\frac{(-1)^{m+1} x^{m+1}}{(m+1)n^m}$$

$$= \frac{1}{2}x\log(2\pi) - \frac{1}{2}x(1+x) - \frac{1}{2}\gamma x^2 - \sum_{m=2}^{\infty}\frac{(-1)^{m+1} x^{m+1}}{m+1}\sum_{n=1}^{\infty}\frac{1}{n^m}$$

$$= \frac{1}{2}x\log(2\pi) - \frac{1}{2}x(1+x) - \frac{1}{2}\gamma x^2 - \sum_{m=2}^{\infty}\frac{(-1)^{m+1}\varsigma(m)x^{m+1}}{m+1}$$

Hence, as reported in [43, p.211], we obtain

(3.17) $$\log G(1+x) = \frac{1}{2}x\log(2\pi) - \frac{1}{2}x(1+x) - \frac{1}{2}\gamma x^2 - \sum_{k=2}^{\infty}\frac{(-1)^{k+1}\varsigma(k)x^{k+1}}{k+1}$$

We now recall the Gosper/Vardi functional equation [6a] for the double gamma function (a further derivation of this is contained in [19])

(3.18) $\quad \log G(1+x) - x\log \Gamma(1+x) = \varsigma'(-1) - \varsigma'(-1, 1+x)$

and using (3.15) and (3.18) we obtain

$$\frac{1}{2}x\log(2\pi) - \frac{1}{2}x(1+x) - \frac{1}{2}\gamma x^2 - \sum_{k=2}^{\infty}\frac{(-1)^{k+1}\varsigma(k)}{k+1}x^{k+1} + \gamma x^2 + \sum_{k=2}^{\infty}\frac{(-1)^{k+1}\varsigma(k)}{k}x^{k+1}$$

$$= \varsigma'(-1) - \varsigma'(-1, 1+x)$$

This may be written as

(3.19) $\quad \varsigma'(-1) - \varsigma'(-1, 1+x) = \frac{1}{2}x\log(2\pi) - \frac{1}{2}x(1+x) + \frac{1}{2}\gamma x^2 + \sum_{k=2}^{\infty}\frac{(-1)^{k+1}\varsigma(k)}{k(k+1)}x^{k+1}$

which is contained in a slightly disguised form in [43, p.222, Eq (539)].

Differentiation results in

$$-\frac{d}{dx}\varsigma'(-1, 1+x) = \frac{1}{2}\log(2\pi) - \frac{1}{2} - x + \gamma x + \sum_{k=2}^{\infty}\frac{(-1)^{k+1}\varsigma(k)}{k}x^k$$

and as shown previously we have

$$\frac{d}{dx}\varsigma'(-1, 1+x) = \varsigma'(0, 1+x) - \varsigma(0, 1+x)$$

$$= \log \Gamma(1+x) - \frac{1}{2}\log(2\pi) + \frac{1}{2} + x$$

Therefore we simply recover (3.15)

$$\log \Gamma(x+1) = -\gamma x + \sum_{k=2}^{\infty}\frac{(-1)^k \varsigma(k)}{k}x^k$$

We have

$$\log S_3(x) = \frac{1}{2}x^2 + \sum_{n=1}^{\infty}\left[n^2 \log\left(1 - \frac{x^2}{n^2}\right) + x^2\right]$$

$$= \frac{1}{2}x^2 + \sum_{n=1}^{\infty}\left[x^2 - \sum_{k=1}^{\infty}\frac{x^{2k}}{kn^{2k-2}}\right]$$

$$= \frac{1}{2}x^2 - \sum_{n=1}^{\infty}\sum_{k=2}^{\infty}\frac{x^{2k}}{kn^{2k-2}} = \frac{1}{2}x^2 - \sum_{k=2}^{\infty}\frac{x^{2k}}{k}\sum_{n=1}^{\infty}\frac{1}{n^{2k-2}}$$

$$= \frac{1}{2}x^2 - \sum_{k=2}^{\infty}\frac{\varsigma(2k-2)}{k}x^{2k}$$

This gives us

(3.21) $$\log S_3(x) = \frac{1}{2}x^2 - \sum_{k=1}^{\infty}\frac{\varsigma(2k)}{k+1}x^{2k+2} = -\sum_{k=0}^{\infty}\frac{\varsigma(2k)}{k+1}x^{2k+2}$$

The following identity for $|x|<1$ appears in the book by Srivastava and Choi [43, p.216, Eq (501)]

$$\sum_{k=1}^{\infty}\frac{\varsigma(2k)}{k+1}x^{2k+2} = \frac{1}{2}[1-\log(2\pi)]x^2 + x\log\frac{G(1+x)}{G(1-x)} - \int_0^x \log G(1+t)\,dt - \int_0^{-x}\log G(1+t)\,dt$$

We note that

$$\int_0^{-x}\log G(1+t)\,dt = -\int_0^x \log G(1-t)\,dt$$

and integration by parts shows that

$$x\log\frac{G(1+x)}{G(1-x)} - \int_0^x \log G(1+t)dt - \int_0^{-x}\log G(1+t)dt = \int_0^x t\frac{d}{dt}\log\frac{G(1+t)}{G(1-t)}dt$$

Using Kinkelin's formula (3.12) we see that

$$\frac{d}{dt}\log\frac{G(1+t)}{G(1-t)} = \log(2\pi) - \pi t \cot \pi t$$

and hence we have

$$\int_0^x t\frac{d}{dt}\log\frac{G(1+t)}{G(1-t)}dt = \frac{1}{2}\log(2\pi)x^2 - \int_0^x \pi t^2 \cot \pi t\, dt$$

This results in

(3.22) $$\sum_{k=1}^{\infty} \frac{\varsigma(2k)}{k+1} x^{2k+2} = \frac{1}{2} x^2 - \int_0^x \pi t^2 \cot \pi t \, dt$$

□

We have the following well-known identity

(3.23) $$t \cot t = \sum_{n=0}^{\infty} (-1)^n \frac{2^{2n} B_{2n}}{(2n)!} t^{2n} \quad , (|t| < \pi)$$

Combining this with Euler's formula

(3.24) $$\varsigma(2n) = (-1)^{n+1} \frac{2^{2n-1} \pi^{2n} B_{2n}}{(2n)!} \quad , (n \geq 1)$$

and, letting $t \to \pi t$, we obtain

(3.25) $$\pi t \cot \pi t = -2 \sum_{n=0}^{\infty} \varsigma(2n) t^{2n} \quad , (|t| < 1)$$

Since the first term of the series (3.23) is equal to $B_0 = 1$, to be consistent with (3.25), we define $\varsigma(0) = -\frac{1}{2}$ (which in fact also coincides with the value determined by the analytic continuation of $\varsigma(s)$ at $s = 0$). We may now multiply (3.25) by $t$ and integrate this to obtain another derivation of (3.22).

We note from [43, p.207] that

$$\int_0^x \log G(1+t) \, dt = \left( \frac{1}{4} - 2 \log A \right) x + \frac{1}{4} \log(2\pi) x^2 - \frac{1}{6} x^3 - (1-x) \log G(1+x) + \log G(x)$$

$$- 2 \log \Gamma_3(1+x) + 2 \log \Gamma_3(x)$$

and hence we may express the triple sine function in terms of the double and triple gamma functions. This type of representation naturally arises from the combination of the facts that

$$\log S_n(x) = \int_0^x \pi t^{n-1} \cot \pi t \, dt$$

$$-\cot \pi t = \psi(t) - \psi(-t) + \frac{1}{t}$$

and the knowledge that integrals of the form $\int_0^x t^{n-1}\psi(t)dt$ result in multiple gamma functions [43, p.208].

The triple gamma function is defined by [17]

$$\Gamma_3(1+x) = \exp(c_1 x + c_2 x^2 + c_3 x^3) F(x)$$

where

$$F(x) = \prod_{k=1}^{\infty} \left\{ \left(1+\frac{x}{k}\right)^{-\frac{1}{2}k(k+1)} \exp\left(\frac{1}{2}(k+1)x - \frac{1}{4}\left(1+\frac{1}{k}\right)x^2 + \frac{1}{6}\frac{x^3}{k} + \frac{1}{6}\frac{x^3}{k^2}\right) \right\}$$

$$= \prod_{k=1}^{\infty} \left\{ \left(1+\frac{x}{k}\right)^{-\frac{1}{2}k(k+1)} \exp\left(\left[1+\frac{1}{k}\right]\left\{\frac{1}{2}kx - \frac{1}{4}x^2 + \frac{1}{6}\frac{x^3}{k}\right\}\right) \right\}$$

and

$$c_1 = \frac{3}{8} - \frac{1}{4}\log(2\pi) - \log A, \qquad c_2 = \frac{1}{4}\left[\gamma + \log(2\pi) + \frac{1}{2}\right]$$

$$c_3 = -\frac{1}{6}\left[\gamma + \varsigma(2) + \frac{3}{2}\right]$$

Journal of Computational and Applied Mathematics 75 (1996) 43-46

Donal F. Connon
Elmhurst
Dundle Road
Matfield, Kent TN12 7HD
dconnon@btopenworld.com